\documentclass[a4paper]{article}
\usepackage[utf8]{inputenc}

\usepackage[english]{babel}
\usepackage{graphicx}
\usepackage{xcolor}
\usepackage{amsmath,amsthm}
\usepackage{relsize}
\usepackage{tikz}
\usetikzlibrary{tikzmark}
\usepackage{amssymb}
\usepackage{mathrsfs}
\usepackage{txfonts,pxfonts,tikz}
\usepackage{tgschola}
\usepackage[T1]{fontenc}
\usepackage[margin= 1 in]{geometry}
\usepackage{mathtools}
\usepackage{enumitem}   
\usepackage{listings}
\usepackage[backend=biber,style=nature,sorting=nyt]{biblatex}
\usepackage{hyperref}
\hypersetup{
    colorlinks=true,
    citecolor=cyan,
    }

\newcommand{\defeq}{\vcentcolon=}

\newcommand{\Q}{\mathbb{Q}}
\newcommand{\C}{\mathbb{C}}
\newcommand{\Z}{\mathbb{Z}}
\newcommand{\N}{\mathbb{N}}
\newcommand{\F}{\mathbb{F}}

\title{1-Bounded Entropy for $C^*$-Algebras}
\author{Connor MacMahon}

\theoremstyle{definition}

\theoremstyle{plain}

\begin{document}

\maketitle

\begin{abstract}
    Towards a question of Voiculescu, two notions of $1$-bounded entropy, $h$ and $h^\text{top}$, are defined for $C^*$-algebras. The quantity $h$ involves the tracial completion with respect to all traces while the quantity $h^\text{top}$ depends on operator norm microstates. It is demonstrated that $h^\text{top}(\mathscr{A})$ does not exceed $h(\mathscr{A})$ for all $C^*$-algebras $\mathscr{A}$. Moreover, a variational principle enabling the computation of $h$ is introduced. The quantity $h$ is computed in various examples, and the computation is used to justify structural conclusions such as $C^*$-primeness, crossed product indecomposability, and free indecomposability. These notions of entropy are generalized to operator systems, where it is demonstrated that both may be computed on a basis.
\end{abstract}

\section{Introduction}

In the 90's Voiculescu introduced the microstates free entropy dimension $\delta_0(x)$ of a tuple $x$ in a tracial von Neumann algebra (\cite{Voic1},\cite{Voic2}). The quantity $\delta_0(x)$ roughly encodes the abundance of matrix tuples which approximately share the same trace statistics with $x$. Following the discovery of $\delta_0$, Voiculescu and others used this quantity to prove striking structural results about von Neumann algebras. For instance, the absence of Cartan subalgebras in free group factors \cite{Voic2}.

\medskip

Despite the utility of $\delta_0$, it is a nagging mystery as to whether $\delta_0$ is dependent on the generating tuple chosen for a given von Neumann algebra. More precisely, if $x$ and $y$ are tuples which generate the same von Neumann algebra, it is unknown whether $\delta_0(x)=\delta_0(y)$. Jung \cite{Jung} introduced the property of strong $1$-boundedness for von Neumann algebras in order to circumvent this issue. Strong $1$-boundedness of a tuple $x$ implies that $\delta_0(x) \leq 1$, and the property of strong $1$-boundedness is an invariant of the generated von Neumann algebra. That is, if $x$ is strongly $1$-bounded and $y$ is such that $W^*(x)=W^*(y)$, then $y$ is also strongly $1$-bounded. 

\medskip

Introduced implicity in the aforementioned work of Jung \cite{Jung} and first written explicitly by Hayes \cite{Reg}, the $1$-bounded entropy $h(M)$ of a tracial von Neumann algebra $M$ has become an invariant of great interest in the past decade. The $1$-bounded entropy is such that $h(M) < \infty$ if and only if $M$ is strongly $1$-bounded. The invariance of $h$ streamlines and expands the applicability of structural results on tracial von Neumann algebras arising from matricial approximations. 

\medskip

Recently, Voiculescu asked whether such an invariant exists in the setting of $C^*$ algebras. In this document, two viable candidates for such a quantity are introduced. Denoted $h$ and $h^\text{top}$ and referred to as $1$-bounded entropy and topological $1$-bounded entropy respectively, the first involves the tracial completion of a $C^*$-algebra with respect to all of its traces, and the second involves the direct consideration of operator norm microstates. This work may be viewed as building on the foundations present in Voiculescu's work \cite{Voic3} in which he developed the proper notions of free entropy and free entropy dimension for $C^*$ algebras. 

\medskip

By definition, both $h$ and $h^\text{top}$ may be computed on a generating set for a $C^*$-algebra by packing microstates spaces. It is demonstated that both quantities are invariants (see Theorem 5.1 and Theorem 7.1). That is, they do not depend on the generating set used to compute them. In comparing the two invariants, one arrives at the following Theorem.

\medskip

\noindent \textbf{Theorem 8.1:} Let $\mathscr{A}$ be a $C^*$-algebra. Then $h^\text{top}(\mathscr{A}) \leq h(\mathscr{A})$. 

\medskip

Moreover, the $1$-bounded entropy $h(\mathscr{A})$ of a $C^*$ algebra $\mathscr{A}$ may be related to the more familiar $1$-bounded entropy of von Neumann algebras via the following result.

\medskip

\noindent \textbf{Theorem 6.1 (Variational Principle):} Let $\mathscr{A}$ be a $C^*$-algebra and $\mathscr{T}(\mathscr{A})$ the collection of all traces on $\mathscr{A}$. For each trace $\tau$ belonging to $\mathscr{T}(\mathscr{A})$, denote by $M_\tau$ the von Neumann algebra generated by the GNS representation of $\mathscr{A}$ corresponding to $\tau$. Then
\begin{align*}
    h(\mathscr{A}) = \sup_{\tau \in \mathscr{T}(\mathscr{A})} h(M_\tau).
\end{align*}

\medskip

This Theorem allows for the computation of the $1$-bounded entropies of $C^*$-algebras in a variety of examples. Moreover, various indecomposability results for $C^*$-algebras are obtained in analogy to the case of tracial von Neumann algebras using this principle.

\medskip

Finally, working in an independent direction, one may define $h$ and $h^\text{top}$ for an operator system $V$ by computing these invariants on the universal $C^*$ algebra $C^*_\text{max}(V)$ of $V$. It is demonstrated that both the $1$-bounded entropy and the topological $1$-bounded entropy of an operator system $V$ may also be computed on a basis for $V$ using microstates appropriate in the category of operator systems, so that both quantities are intrinsic to the operator system $V$. 

\medskip

\noindent \textbf{Acknowledgement:} The author acknowledges support from NSF Career grant DMS-21447. Additionally, the author wishes to recognize the unwavering and substantial support granted by his advisor Ben Hayes, without which this project would not exist. The author also thanks David Jekel for helpful feedback as version 2 was written. 

\section{Tracially Complete $C^*$-Algebras}

\textbf{Definition 2.1:} Fix a unital $C^*$-algebra $M$ and a weak-$*$ compact, convex collection of traces $\mathscr{T}$ on $M$. The pair $(M,\mathscr{T})$ is said to be a tracially complete $C^*$-algebra provided the expression 
\begin{align*}
    \|a\|_{2,\mathscr{T}} \defeq \sup_{\tau \in \mathscr{T}}\tau(a^*a)^{1/2}
\end{align*}
defines a norm on $M$ such that $\text{Ball}(M)$ is complete in the metric induced by $\|\cdot\|_{2,\mathscr{T}}$ \cite{TC}.

\medskip

\noindent \textbf{Remark:} That $\|\cdot\|_{2,\mathscr{T}}$ defines a norm on $M$ amounts to the faithfulness of:
\begin{align*}
    \bigoplus_{\tau \in \mathscr{T}} \{\pi_\tau,L^2(M,\tau)\},
\end{align*}
the direct sum of GNS representations over all traces in $\mathscr{T}$.

\medskip

Observe that in the case of a single trace $\tau$, the pair $(M,\tau)$ is a tracial von Neumann algebra. Hence, these objects are viewed as generalizing von Neumann algebras. 

\medskip

The following analogue of the Kaplansky Density Theorem is both useful and illustrates the analogy quite well. Henceforth, this result is referred to simply as the Kaplansky Density Theorem.

\medskip

\noindent \textbf{Theorem (\cite{TC} 3.28):} Let $(M,\mathscr{T})$ be a tracially complete $C^*$-algebra and let $\mathscr{A}$ be a $\|\cdot\|_{2,\mathscr{T}}$-dense $C^*$-subalgebra of $M$. Then the unit ball of $\mathscr{A}$ is $\|\cdot\|_{2,\mathscr{T}}$-dense in the unit ball of $M$. 

\medskip

Now, the primary mechanism by which tracially complete $C^*$ algebras arise in practice is presented. This construction was initially carried out in work of Ozawa \cite{Ozawa} and elaborated upon by Carrión et al. \cite{TC}. 

\medskip

\noindent \textbf{Example 2.1:} Let $\mathscr{A}$ be a $C^*$-algebra and $\mathscr{T}$ a weak-$*$ compact, convex collection of traces on $\mathscr{A}$. Define the tracial completion $\overline{\mathscr{A}}^\mathscr{T}$ as follows:
\begin{align*}
    \overline{\mathscr{A}}^\mathscr{T} = \frac{\{(a_n)_{n \in \N} \in \ell^\infty(\mathscr{A},\N) \, | \, (a_n)_{n \in \N} \, \text{is} \, \|\cdot\|_{2,\mathscr{T}} \, \text{Cauchy}\}}{\{(a_n)_{n \in \N} \in \ell^\infty(\mathscr{A},\N) \, | \, \|a_n\|_{2,\mathscr{T}} \rightarrow 0\}}.
\end{align*}
Completeness of the ball of $\overline{\mathscr{A}}^\mathscr{T}$ follows essentially by construction. Moreover, each trace in $\mathscr{T}$ lifts to a trace on $\overline{\mathscr{A}}^\mathscr{T}$ as applying a trace to a member of $\overline{\mathscr{A}}^\mathscr{T}$ term by term yields a Cauchy sequence in $\C$. Simply declare the output of the trace to be the limit. Finally, $\mathscr{A}$ embeds diagonally within $\overline{\mathscr{A}}^\mathscr{T}$. The tracially complete $C^*$ algebra $\overline{\mathscr{A}}^\mathscr{T}$ is called the tracial completion of the pair $(\mathscr{A},\mathscr{T})$. When $\mathscr{T}$ is taken to be the collection of all traces on $\mathscr{A}$, then $\overline{\mathscr{A}}^\mathscr{T}$ is referred to as the universal tracial completion of $\mathscr{A}$.

\medskip

In order to define an analogue of $1$-bounded entropy for $C^*$-algebras, the method employed below is to pass first to the universal tracial completion. 

\medskip

It is proper to warn the reader that this approach only addresses tracial $C^*$ algebras. For instance, the Calkin algebra possesses no traces, where both invariants $h$ and $h^\text{top}$ evaluate to $-\infty$. In the following considerations involving $1$-bounded entropy it will be tacitly assumed that the $C^*$ algebras in question admit traces. 

\section{The Space of Laws}

Before entropy is defined, the notion of a law must be introduced. Most of the subsequent exposition on spaces of laws and entropy follows the treatment in Hayes, Jekel, Nelson, and Sinclair \cite{Absorption}.  

\medskip

\noindent \textbf{Definition 3.1 (\cite{Absorption} section 2.1):} Let $I$ be an index set. A functional $\ell: \C^*\langle T_i \, | i \in I \rangle \rightarrow \C$ on the polynomial algebra on $I$ free variables and their formal adjoints is called a tracial law provided there is $R>0$ such that:
\begin{enumerate}
    \item $\ell(P^*P) \geq 0$ for all members $P$ of $\C^*\langle T_i \, | \, i \in I \rangle$,
    \item $\ell(1)=1$,
    \item $\ell(PQ)=\ell(QP)$ for all members $P,Q$ of $\C^*\langle T_i \, | \, i \in I \rangle$,
    \item and for each natural number $k$ and member $j$ of $I$ there holds:
    \begin{align*}
        \ell([T_j^*T_j]^{k}) \leq R^{2k}.
    \end{align*}
\end{enumerate}
Denote by $\Sigma_I$ the space of such tracial laws. Given a constant $R>0$, set $\Sigma_{R,I}$ to be the collection of all tracial laws $\ell$ satisfying (4) for the constant $R$.

\medskip

Since the space of laws is a subset of the algebraic dual of $\C^*\langle T_i \, | \, i \in I \rangle$, it may be equipped with the weak-$*$ topology. In this topology, $\Sigma_{R,I}$ is compact and Hausdorff. This is readily seen via an embedding of $\Sigma_{R,I}$ into a Tychonoff product indexed by $*$-polynomials.

\medskip

Each tuple $x \defeq (x_i)_{i \in I}$ within a tracial von Neumann algebra $(M,\tau)$ with $\|x_i\| \leq R$ for all $i$ defines a member $\ell_x$ of $\Sigma_{R,I}$ as follows:
\begin{align*}
    \ell_x(P) = \tau(P(x)).
\end{align*}
Conversely, each member of $\Sigma_{R,I}$ arises in this way by a variant of the GNS construction. Hence, $\Sigma_{R,I}$ parametrizes $R$-bounded $I$-tuples in tracial von Neumann algebras up to trace statistics. 

\medskip

Now, suppose that one fixes a generating tuple $a=(a_i)_{i \in I}$ within a tracially complete $C^*$-algebra $(M,\mathscr{T})$ (I.e. the $C^*$-subalgebra generated by $a$ is $\|\cdot\|_{2,\mathscr{T}}$-dense in $M$) which is $R$-bounded. Each trace $\tau$ in $\mathscr{T}$ defines a law as follows:
\begin{align*}
    \ell_\tau(P) = \tau(P(a)). 
\end{align*}
In this way, the set $\mathscr{T}$ gives rise to the following compact, convex subset of $\Sigma_{R,I}$:
\begin{align*}
    L_a \defeq \{\ell \in \Sigma_{R,I} \, | \, \exists \tau \in \mathscr{T} \, \text{such that} \, \ell(P) = \tau(P(a))\}. 
\end{align*}
In the tracial von Neumann algebra setting, the entropy of a tuple itself is calculated. The above necessitates that the entropy of a compact, convex subset of the space of laws be defined. This is the objective of the following section. 

\section{Definition of The Entropy}

Let $K$ be a compact, convex subset of the space $\Sigma_{R,I}$ of $R$-bounded laws. Fix a natural number $N$ and a weak-$*$ neighborhood $\mathcal{O}$ of $K$ in $\Sigma_{R,I}$. Given these inputs, one defines the associated microstate space $\Gamma_{R}^{(N)}(\mathcal{O})$ as follows (\cite{Absorption} section 2.2):
\begin{align*}
    \Gamma_{R}^{(N)}(\mathcal{O}) \defeq \{A=(A_i)_{i \in I} \in M_N(\C)^I \, | \, \|A_i\| \leq R \, \text{ for all $i$, and} \, \ell_A \in \mathcal{O} \}.
\end{align*}
Observe that the microstate space is stable under unitary conjugation. In order to define the entropy as in the work of Hayes, Jekel, Nelson, and Sinclair \cite{Absorption} (see Definition 2.2), one needs to pack these microstates spaces in $2$-norm. The $2$-norm on $M_N(\C)$ is defined as follows:
\begin{align*}
    \|A\|_2 \defeq [\frac{1}{N}\text{Tr}(A^*A)]^{1/2}. 
\end{align*}
Make a note that the normalized trace is used. Now, fix a subset $\Omega$ of $M_N(\C)^I$. For each tuple $B$ in $\Omega$, each finite subset $F \subset I$, and each $\epsilon>0$, one defines the $(F,\epsilon)$ orbital neighborhood $N_{F,\epsilon}^\text{orb}(B)$ of $B$ as follows (\cite{Absorption} 2.2):
\begin{align*}
    N_{F,\epsilon}^\text{orb}(B) \defeq \{C \in M_N(\C)^I \, | \, \exists U \in U(N) \, \text{such that} \, \|C_i-UB_iU^*\|_2 < \epsilon \, \forall i \in F\},
\end{align*}
where $U(N) \subset M_N(\C)$ denotes the unitary group. If no "orb" superscript is specified, it is simply the collection of $C$ with $\|C_i-B_i\|_2<\epsilon$. For a collection $\Omega_0 \subset \Omega$ of members of $\Omega$, denote (\cite{Absorption} 2.2):
\begin{align*}
    N_{F,\epsilon}^\text{orb}(\Omega_0) \defeq \bigcup_{B \in \Omega_0} N_{F,\epsilon}^\text{orb}(B).
\end{align*}
A subset $\Omega_0$ of $\Omega$ is said to be $(F,\epsilon)$ orbit dense provided $\Omega \subset N_{F,\epsilon}^\text{orb}(\Omega_0)$. Now, one may define the $(F,\epsilon)$ orbital covering number $K_{F,\epsilon}^\text{orb}(\Omega)$ of $\Omega$ to be the minimal cardinality of an $(F,\epsilon)$ orbit dense subset of $\Omega$. That is:
\begin{align*}
    K_{F,\epsilon}^\text{orb}(\Omega) \defeq \inf \{|\Omega_0| \, | \, \Omega_0 \subset \Omega \, \text{and} \, \Omega \subset N_{F,\epsilon}^\text{orb}(\Omega_0)\}. 
\end{align*}
By compactness, this quantity will be finite provided $\Omega$ is contained within a product of $\|\cdot\|_2$-balls. 

\medskip

Now, it's time to define the entropy. The definition is quite similar to that in Hayes, Jekel, Nelson, and Sinclair \cite{Absorption}, with one notable difference: instead of taking weak-$*$ neighborhoods of a single law, weak-$*$ neighborhoods of an entire compact, convex set $K$ of laws are considered. 

\medskip

\noindent \textbf{Definition 4.1:} Let $K$ be a compact, convex subset of $\Sigma_{R,I}$. The $1$-bounded entropy of $K$, denoted $h(K)$, is defined as follows:
\begin{align*}
    h(K) = \sup_{F,\epsilon} \inf_{\mathcal{O} \supset K} \limsup_{N \rightarrow \infty} \frac{1}{N^2}\log (K_{F,\epsilon}^\text{orb}(\Gamma_R^{(N)}(\mathcal{O}))).
\end{align*}
Here, $F$ denotes a finite subset of $I$ and $\mathcal{O}$ a weak-$*$ neighborhood of $K$ in $\Sigma_{R,I}$. 

\medskip

Similarly, one may define the $1$-bounded entropy of an inclusion of tracially complete $C^*$ algebras $\mathscr{A} \subset \mathscr{B}$.

\medskip

\noindent \textbf{Definition 4.2:} Suppose that $\mathscr{A} \subset \mathscr{B}$ is an inclusion of tracially complete $C^*$-algebras. Fix an $R$-bounded generating tuple $(a_i)_{i \in I}$ for $\mathscr{A}$ and an $R$-bounded tuple $(b_j)_{j \in J}$ of elements in $\mathscr{B}$ so that $\mathscr{B}$ is jointly generated by $(a_i)_{i \in I}$ and $(b_j)_{j \in J}$. Consider the compact, convex set $L_{(a,b)}$ of laws within $\Sigma_{R,I \sqcup J}$, and define:
\begin{align*}
    h(\mathscr{A}:\mathscr{B}) \defeq \sup_{F \subset I,\epsilon} \inf_{\mathcal{O} \supset L_{a,b}} \limsup_{N \rightarrow \infty} \frac{1}{N^2}\log (K_{F,\epsilon}^\text{orb}(\Gamma_R^{(N)}(\mathcal{O}))).
\end{align*}
Here, in the supremum $F$ ranges over finite subsets of the index set $I$. The sets $\mathcal{O}$ are open neighborhoods of $L_{(a,b)}$ within $\Sigma_{R,I \sqcup J}$. This quantity is referred to as the $1$-bounded entropy of $\mathscr{A}$ in the presence of $\mathscr{B}$. Observe that $h(L_a)$ as defined above is the $1$-bounded entropy $h(\mathscr{A}:\mathscr{A})$ of $\mathscr{A}$ in the presence of itself. This quantity is referred to simply as the $1$-bounded entropy of $\mathscr{A}$ and denoted $h(\mathscr{A})$. 

\medskip

The $1$-bounded entropy $h(\mathscr{A}: \mathscr{B})$ of $\mathscr{A}$ in the presence of $\mathscr{B}$ is not yet well defined: a priori, the quantity may depend on the choice of the tuples $a$ and $b$. Proving that $h(\mathscr{A}: \mathscr{B})$ is in fact well defined is the business of the next section. Until this proof is done, the quantity will be denoted $h(a:b)$, specifying the tuples used in its computation. Observe that in the case of a tracial von Neumann algebras, the definition reduces to the definition of the $1$-bounded entropy in the presence given in work of Hayes \cite{Reg} (section 2). 

\section{Proof of Invariance}

The objective of this section is to prove the following result, which is the analogue of the invariance claim in Hayes' work \cite{Reg} (appendix A).

\medskip

\noindent \textbf{Theorem 5.1:} Suppose that $\mathscr{A} \subset \mathscr{B}$ is an inclusion of tracially complete $C^*$-algebras. Fix positive constants $R,S>0$ and index sets $I_1,J_1$ and $I_2,J_2$. Suppose $a\defeq(a_i)_{i \in I_1}$ and $b\defeq(b_i)_{i \in J_1}$ are $R$-bounded such that $a$ generates $\mathscr{A}$ and $a$ and $b$ jointly generate $\mathscr{B}$. Suppose similarly that $c=(c_i)_{i \in I_2}$ and $d \defeq (d_i)_{i \in J_2}$ are $S$-bounded so that $c$ generates $\mathscr{A}$ and $c$ and $d$ jointly generate $\mathscr{B}$. Then
\begin{align*}
    h(a:b) = h(c:d).
\end{align*}
That is, the $1$-bounded entropy $h(\mathscr{A}:\mathscr{B})$ of $\mathscr{A}$ in the presence of $\mathscr{B}$ is independent of the tuples used in its computation. 

\medskip

In building up to this result, first a definition is required. It will be useful to have a universal $C^*$ algebra whose quotients parametrize $C^*$-algebras generated by $R$-bounded tuples. More precisely, the following definition is what's required (see Hayes, Jekel, Nelson, and Sinclair \cite{Absorption} section 3 for a similar universal construction).

\medskip

\noindent \textbf{Definition 5.1:} Fix an index set $I$ and $R>0$. Consider the product:
\begin{align*}
    \Pi \defeq \prod_{i \in I} R\text{Ball}(B(\ell^2(I))).
\end{align*}
If $P$ is a member of $\C^*\langle T_i \, | \, i \in I\rangle$, define:
\begin{align*}
    \|P\|_R \defeq \sup_{x \in \Pi} \|P(x)\| < \infty.
\end{align*}
Here, $P(x)$ evidently denotes the evaluation of $P$ on the tuple $x$, substituting $x_i$ for $T_i$. Observe that $\|\cdot\|_R$ defines a seminorm on $\C^*\langle T_i \, | \, i \in I\rangle$. Denote by $C^*(R)$ the separation and completion of $\C^*\langle T_i \, | \, i \in I\rangle$ with respect to $\|\cdot\|_R$, and observe that $C^*(R)$ has the following universal property: for each unital $C^*$ algebra $\mathscr{A}$ and each $R$-bounded $I$ tuple $x = (x_i)_{i \in I}$ of members of $\mathscr{A}$, the evaluation $P \mapsto P(x)$ extends uniquely to a *-homomorphism $C^*(R) \rightarrow \mathscr{A}$. 

\medskip

This definition enables the formulation and use of norm estimates on evaluations which hold uniformly over all $R$-bounded $I$ tuples. For example, the following Corollary.

\medskip

\noindent \textbf{Corollary 5.1:} Let $(M,\mathscr{T})$ be a tracially complete $C^*$-algebra. Let $I$ and $J$ be index sets, $R,S>0$, and suppose that $a \in M^I$ and $b \in M^J$ are $R$ and $S$ bounded respectively. Then for each $\epsilon>0$ and $F \subset J$ finite, there exists an $F$-tuple $(P_j)_{j \in F}$ of polynomials in $\C^*\langle T_i \, | \, i \in I \rangle$ such that:
\begin{align*}
    \|P_j\|_R \leq S
\end{align*}
for each $j \in F$, and
\begin{align*}
    \|P_j(a)-b_j\|_{2,\mathscr{T}} < \epsilon
\end{align*}
for each $j \in F$. 
\begin{proof}
    Suppose without loss of generality that $S=1$. Denote by $\phi$ the $*$-homomorphism $C^*(R) \rightarrow M$ extending the evaluation on $a$. Since $a$ generates $M$, one has that $\phi(C^*(R))$ is a $\|\cdot\|_{2,\mathscr{T}}$-dense $C^*$-subalgebra of $M$. Hence, the Kaplansky density theorem \cite{TC} yields the following:
    \begin{align*}
        \overline{\text{Ball}(\phi(C^*(R)))}^{\|\cdot\|_{2,\mathscr{T}}} = \text{Ball}(M). 
    \end{align*}
    Since $\phi$ induces an isometry $C^*(R)/\ker \phi \rightarrow M$ on the factor space $C^*(R)/\ker \phi$, one has that $\text{Ball}^0(\phi(C^*(R))) = \phi(\text{Ball}^0(C^*(R)))$. Recall that the "$0$" superscript denotes the interior. Therefore, since polynomials form a norm-dense $*$-subalgebra of $C^*(R)$, one may choose polynomials of norm at most $1$ such that their evaluations on $a_j$ for $j \in F$ approximate $b_j$ within $\epsilon$ in the $\|\cdot\|_{2,\mathscr{T}}$ norm. 
\end{proof}

\medskip

The above corollary gives a nice means of transferring between norm bounded tuples. However, one must control neighborhoods of the set $L_a$ under such a change. The following Lemma accounts for this.

\medskip

\noindent \textbf{Remark:} In the following proof, it is used tacitly that the $\|\cdot\|_{2,E}$-$\|\cdot\|_2$ Lipschitz constant of a $*$-polynomial on $(R\text{Ball}(M_N(\C)))^I$ is independent of $N$ for each finite subset $E$ of $I$. In fact, there is a universal Lipschitz constant which applies over all $R$-bounded $I$ tuples in arbitrary tracial von Neumann algebras. Indeed, suppose for contradiction that $(M_n,\tau_n)$ is a sequence of tracial von Neumann algebras with the supremum of their associated Lipschitz constants infinity, then their free product
\begin{align*}
     \mathlarger{\mathlarger{\mathlarger{\ast}}}_{n \in \N} (M_n,\tau_n)
\end{align*}
produces a single tracial von Neumann algebra for which $P$ is not $\|\cdot\|_{2,E}$-$\|\cdot\|_2$ Lipschitz on $R$-bounded $I$ tuples. The Lipschitz claim on a single tracial von Neumann algebra amounts to a reduction to monomials via the triangle inequality and then an induction on degree (see Hayes, Jekel, and Kunnawalkam Elayavalli \cite{Betti} section 4 for the self-adjoint case).

\medskip

\noindent \textbf{Lemma 5.1:} Let $(M,\mathscr{T})$ be a tracially complete $C^*$-algebra. Let $I$ and $J$ be index sets, $R,S>0$, and suppose that $a \in M^I$ is an $R$-bounded tuple. Recall that:
\begin{align*}
   & L_a \defeq \{\ell \in \Sigma_{R,I} \, | \, \exists \tau \in \mathscr{T} \, \text{such that} \, \ell(P) = \tau(P(a))\}.
\end{align*}
Suppose that $Q = (Q_j)_{j \in J}$ is a tuple of polynomials in $\C^*\langle T_i \, | \, i \in I \rangle$ such that $\|Q_j\|_R \leq S$ for each $j$. Then each of the following holds.
\begin{enumerate}
    \item For each weak-$*$ neighborhood $\mathcal{O}$ of $L_{Q(a)}$ there exists a weak-$*$ neighborhood $\mathcal{V}$ of $L_a$ such that:
    \begin{align*}
        Q(\Gamma_R^{(N)}(\mathcal{V})) \subset \Gamma_S^{(N)}(\mathcal{O})
    \end{align*}
    for all $N$.
    \item For each weak-$*$ neighborhood $\mathcal{O}$ of $L_a$ and finite $F \subset I$, there exists a finite subset $E$ of $I$ containing $F$, a constant $\delta >0$, and a weak-$*$ neighborhood $\mathcal{V}$ of $L_a$ such that:
    \begin{align*}
        N_{E,\delta}(\Gamma_R^{(N)}(\mathcal{V})) \cap (R\text{Ball}(M_N(\C)))^I \subset \Gamma_R^{(N)}(\mathcal{O}). 
    \end{align*}
    \item For each weak-$*$ neighborhood $\mathcal{O}$ of $L_a$, there is a finite index set $F \subset I$ and a constant $\eta >0$ such that if $\tilde{a}$ is $R$-bounded with 
    \begin{align*}
        \|a_j-\tilde{a}_j\|_{2,\mathscr{T}} < \eta,
    \end{align*}
    for all $j$ belonging to $F$, then there is a weak-$*$ neighborhood $\mathcal{V}$ of $L_{\tilde{a}}$ satisfying:
    \begin{align*}
        \Gamma_R^{(N)}(\mathcal{V}) \subset \Gamma_R^{(N)}(\mathcal{O})
    \end{align*}
    for all $N$. 
\end{enumerate}
\begin{proof}
    First, a proof of (1). Let $E$ be a finite collection of polynomials in $\C^*\langle T_j \, | \, j \in J\rangle$ and fix $\epsilon > 0$. Denote by $U_{E,\epsilon}$ the weak-$*$ sub-basic neighborhood defined by:
    \begin{align*}
        U_{E,\epsilon} \defeq \{\ell \in \C \langle T_j \, | \, j \in J \rangle' \, | \, |\ell(P)|<\epsilon \, \text{for all} \, P \in E\}. 
    \end{align*}
   Here, $\C \langle T_j \, | \, j \in J \rangle'$ denotes the algebraic dual of $\C \langle T_j \, | \, j \in J \rangle$. By compactness of $L_a$, without loss of generality $\mathcal{O} = (L_{Q(a)} + U_{E,\epsilon}) \cap \Sigma_{R,I}$. Now, set
    \begin{align*}
        \mathcal{V} \defeq \bigcap_{P \in E} \bigcup_{\tau \in \mathscr{T}} \{\ell \in \Sigma_{R,I} \, | \, |\ell(P \circ Q) - \tau(P(Q(a)))|<\epsilon\}.
    \end{align*}
    Evidently $\mathcal{V}$ is weak-$*$ open, and contains $L_a$. Now, recall that a tuple $A =(A_i)_{i \in I}$ of $N \times N$ matrices belongs to $\Gamma_R^{(N)}(\mathcal{V})$ provided it is $R$-bounded and the law $\ell_A$ belongs to $\mathcal{V}$. Since $\|Q_j\|_R \leq S$ for all $j$, the image $Q(A)$ is $S$-bounded. Moreover, since $\ell_A$ belongs to $\mathcal{V}$ there exists a member $\tau$ of $\mathscr{T}$ such that:
    \begin{align*}
        |\ell_A(P \circ Q) - \tau(P(Q(a)))| < \epsilon
    \end{align*}
    for each $P$ in $E$. But this in turn implies that:
    \begin{align*}
        |\ell_{Q(A)}(P) - \tau(P(Q(a)))|< \epsilon.
    \end{align*}
    Hence, $\ell_{Q(A)}$ belongs to $\mathcal{O}$ by definition. 
\end{proof}

\medskip

Now, for a proof of (2). First, some notation. For a finite subset $G$ of $I$, denote by $\|\cdot\|_{2,G}$ the following seminorm on matrix tuples:
\begin{align*}
    \|(B_i)_{i \in I}\|_{2,G} \defeq \bigg( \sum_{i \in G} \|B_i\|_2^2 \bigg)^{1/2}.
\end{align*}
The $2$-norm on individual matrices is that coming from the normalized trace. Once again appealing to compactness, one may set $\mathcal{O} = (L_a + U_{Q,\epsilon}) \cap \Sigma_{R,I}$ for $U_{Q,\epsilon}$ a weak-$*$ sub-basic open neighborhood as above with $Q$ a finite set of $*$-polynomials. Observe that each member of $Q$ is a linear combination of finitely many monomials in the variables $T_j$ and their adjoints. Adjoin to $F$ all indices $j$ for which either $T_j$ or $T_j^*$ appears within a monomial in some member of $Q$, and call the new set $E$. 

\medskip

Set $\mathcal{V} \defeq L_a + U_{Q,\epsilon/2}$. Suppose now that $\delta>0$ is fixed, and choose an $R$-bounded tuple $B \defeq (B_i)_{i \in I}$ of $N \times N$ matrices belonging to $N_{E,\delta}(\Gamma_R^{(N)}(\mathcal{V}))$. By definition, there exists a tuple $A$ in $\Gamma_R^{(N)}(\mathcal{V})$ such that $\|A_i-B_i\|_2<\delta$ for each $i \in E$. Since $E$ is finite, each member $P$ of $Q$ is $\|\cdot\|_{2,E}$-$\|\cdot\|_2$ Lipschitz on $(R\text{Ball}(M_N(\C)))^I$. Denote the Lipschitz constant $\|P\|_{\text{Lip}}$ (see the remark preceding the Lemma). Since $E$ includes all indices corresponding to variables which appear nontrivially in members of $Q$, one has:
\begin{align*}
    \bigg|\frac{1}{N}\text{Tr}(P(A))-\frac{1}{N}\text{Tr}(P(B))\bigg| \leq \|P(A)-P(B)\|_{2} \leq \|P\|_{\text{Lip}} \cdot \|A-B\|_{2,E} \leq \|P\|_{\text{Lip}} \cdot |E|^{1/2} \cdot \delta. 
\end{align*}
Hence, if
\begin{align*}
    \delta = \min_{P \in Q} \epsilon(2\|P\|_{\text{Lip}} \cdot |E|^{1/2} +1)^{-1},
\end{align*}
one has for each member $P$ of $Q$ that:
\begin{align*}
    \bigg|\frac{1}{N}\text{Tr}(P(A))-\frac{1}{N}\text{Tr}(P(B))\bigg|< \epsilon/2.
\end{align*}
On the other hand, since $A$ is a member of $\Gamma_R^{(N)}(\mathcal{V})$, for each member $P$ of $Q$ there exists a trace $\tau$ belonging to $\mathscr{T}$ such that 
\begin{align*}
    \bigg|\frac{1}{N}\text{Tr}(P(A))-\tau(P(a))\bigg| < \epsilon/2.
\end{align*}
Since $P$ was arbitrary, applying the triangle inequality yields that $B$ is a member of $\Gamma_R^{(N)}(\mathcal{O})$. 

\medskip

To demonstrate (3), once again appeal to compactness to take $\mathcal{O} = (L_a + U_{Q,\epsilon}) \cap \Sigma_{R,I}$ where $Q$ is a finite collection of $*$-polynomials. Similarly to (2), consider the following seminorm for each finite subset $G$ of $I$:
\begin{align*}
    \|(b_i)_{i \in I}\|_{2,G,\mathscr{T}} \defeq \bigg(\sum_{i \in G} \|b_i\|_{2,\mathscr{T}}^2 \bigg)^{1/2}. 
\end{align*}
Take $F$ to be the collection of indices $j$ such that either $T_j$ or $T_j^*$ appears nontrivially in some monomial summand of a member of $Q$. Since $F$ is finite, one has that each member $P$ of $Q$ is $\|\cdot\|_{2,F,\mathscr{T}}$-$\|\cdot\|_{2,\mathscr{T}}$ Lipschitz on $(R\text{Ball}(M))^I$. Denote the Lipschitz constant $\|P\|_{\text{Lip}}$. Set $\mathcal{V} \defeq (L_{\Tilde{a}} + U_{Q,\epsilon/2}) \cap \Sigma_{R,I}$. Suppose that $\eta>0$ has been fixed and that $\Tilde{a}$ is $R$-bounded such that:
\begin{align*}
    \|a_j-\Tilde{a}_j\|_{2,\mathscr{T}} < \eta
\end{align*}
for each index $j$ in $F$. Now, for each member $P$ of $Q$ and trace $\tau$ belonging to $\mathscr{T}$ one has the following estimate by Cauchy-Schwarz:
\begin{align*}
    |\tau(P(a))-\tau(P(\Tilde{a}))| \leq \|P(a)-P(\Tilde{a})\|_{2,\tau} \leq \|P(a)-P(\Tilde{a})\|_{2,\mathscr{T}} \leq \|P\|_{\text{Lip}} \cdot \|a-\Tilde{a}\|_{2,F,\mathscr{T}} \leq \|P\|_\text{Lip} \cdot |F|^{1/2} \cdot \eta.
\end{align*}
Hence, if
\begin{align*}
    \eta = \min_{P \in Q} \epsilon (2\|P\|_\text{Lip} \cdot |F|^{1/2} + 1)^{-1},
\end{align*}
then for each member $P$ of $Q$ one has:
\begin{align*}
     |\tau(P(a))-\tau(P(\Tilde{a}))| < \epsilon/2.
\end{align*}
Now, suppose that $A$ belongs to $\Gamma_R^{(N)}(\mathcal{V})$. Then there exists a trace $\tau$ in $\mathscr{T}$ such that:
\begin{align*}
    \bigg|\frac{1}{N} \text{Tr}(P(A)) - \tau(P(\Tilde{a}))\bigg| < \epsilon/2
\end{align*}
for each member $P$ of $Q$. Applying the triangle inequality, one sees that $A$ also belongs to $\Gamma_R^{(N)}(\mathcal{O})$. This completes the proof, with the above choice of $\eta>0$. 

\medskip

Now, the invariance claim may be proved. The reader may wish to consult again the definition of the entropy before reading the proof.

\medskip

\noindent \textbf{Theorem 5.1:} Suppose that $\mathscr{A} \subset \mathscr{B}$ is an inclusion of tracially complete $C^*$-algebras. Fix positive constants $R,S>0$ and index sets $I_1,J_1$ and $I_2,J_2$. Suppose $a\defeq(a_i)_{i \in I_1}$ and $b\defeq(b_i)_{i \in J_1}$ are $R$-bounded such that $a$ generates $\mathscr{A}$ and $a$ and $b$ jointly generate $\mathscr{B}$. Suppose similarly that $c=(c_i)_{i \in I_2}$ and $d \defeq (d_i)_{i \in J_2}$ are $S$-bounded so that $c$ generates $\mathscr{A}$ and $c$ and $d$ jointly generate $\mathscr{B}$. Then
\begin{align*}
    h(a:b) = h(c:d).
\end{align*}
That is, the $1$-bounded entropy $h(\mathscr{A}:\mathscr{B})$ of $\mathscr{A}$ in the presence of $\mathscr{B}$ is independent of the tuples used in its computation. 
\begin{proof}
Fix a finite subset $F$ of the index set $I_1$ and $\epsilon>0$. Applying the Corollary to $(\mathscr{B},\mathscr{T})$, the index sets $I_1,I_2$, and norm cutoffs $R,S>0$, there exists a tuple $(P_i)_{i \in F}$ of $*$-polynomials in $\C^* \langle T_j \, | \, j \in I_2 \rangle$ such that:
\begin{align*}
    \|P_i(c)-a_i\|_{2,\mathscr{T}} < \epsilon.
\end{align*}
for each $i$ belonging to $F$. Additionally, the $P_i$ may be taken to satisfy $\|P_i\|_S \leq R$ for each $i$ in $F$. Now, observe that each $P_i$ is $\|\cdot\|_{2,H,\mathscr{T}}$-$\|\cdot\|_{2,\mathscr{T}}$ Lipschitz on $(S\text{Ball}(\mathscr{B}))^{I_2}$ for each finite subset $H$ of $I_2$ containing all indices $j$ of variables $T_j$ or $T_j^*$ appearing nontrivially in the $P_i$. Recall that for an $I_2$-tuple $x$:
\begin{align*}
    \|x\|_{2,H,\mathscr{T}} \defeq \bigg( \sum_{i \in H} \|x_i\|_{2,\mathscr{T}} \bigg)^{1/2}.
\end{align*}
Let $L'(H)$ denote the maximum of the Lipschitz constants among the $P_i$, so that for each $i$ belonging to $F$ there holds:
\begin{align*}
    \|P_i(x)-P_i(y)\|_{2,\mathscr{T}} \leq L'(H) \cdot \| x-y \|_{2,H,\mathscr{T}}
\end{align*}
for all $S$-bounded $I_2$ tuples $x$ and $y$ in $\mathscr{B}$ and $H$ as above. For convenience, adjoin to $(P_i)_{i \in F}$ additional terms to make the evaluation on $c$ into an $R$-bounded $I_1 \sqcup J_1$ tuple, denoted $P(c)$. For this purpose, one may take the $P_i$ for $i$ not in $F$ to be $0$. 

\medskip

Now, fix a weak-$*$ neighborhood $\mathcal{O}$ of $L_{(c,d)}$ within $\Sigma_{S,I_2 \sqcup J_2}$. Applying (3) in the above Lemma, there exists $\eta >0$ and a finite subset $E$ of $I_2 \sqcup J_2$ so that if $(\Tilde{c},\Tilde{d})$ is an $S$-bounded $I_2 \sqcup J_2$ tuple satisfying 
\begin{align*}
    \|(c,d)_i - (\Tilde{c},\Tilde{d})_i\|_{2,\mathscr{T}} < \eta
\end{align*}
for each $i$ in $E$, then there is a weak-$*$ neighborhood $\mathcal{V}$ of $L_{(\Tilde{c},\Tilde{d})}$ so that:
\begin{align*}
    \Gamma_S^{(N)}(\mathcal{V}) \subset \Gamma_S^{(N)}(\mathcal{O})
\end{align*}
for all $N$. Without loss of generality, suppose that $E$ involves all indices $j$ from $I_2$ corresponding to variables $T_j$ or $T_j^*$ which appear nontrivially in the $P_i$. Also without loss of generality, suppose $\eta < \epsilon/[|E \cap I_2|^{1/2}(L'(E \cap I_2) + 1)]$.

\medskip

Applying the Corollary to $(\mathscr{B},\mathscr{T})$, the index sets $I_2 \sqcup J_2,I_1 \sqcup J_1$, and norm cutoffs $S,R>0$, there exists a tuple $(Q_j)_{j \in E}$ of $*$-polynomials in $\C^* \langle T_i \, | \, i \in I_1 \sqcup J_1 \rangle$ such that:
\begin{align*}
    \|Q_j(a,b) - (c,d)_j\|_{2,\mathscr{T}} < \eta
\end{align*}
for each $j$ belonging to $E$. Moreover, the $Q_j$ may be taken so that $\|Q_j\|_R \leq S$. By adding additional $0$ polynomials $Q_j$ for $j$ not in $E$, complete $(Q_j(a,b))_{j \in E}$ to an $S$-bounded $I_2 \sqcup J_2$ tuple denoted $Q(a,b)$. By the definition of $\eta$, there is a weak-$*$ neighborhood $\mathcal{V}$ of $L_{Q(a,b)}$ satisfying:
\begin{align*}
    \Gamma_S^{(N)}(\mathcal{V}) \subset \Gamma_S^{(N)}(\mathcal{O})
\end{align*}
for all $N$. Observe the following for each $i$ in $F$:
\begin{align*}
    \|P_i(Q(a,b))-a_i\|_{2,\mathscr{T}} = \|P_i(Q(a,b)) - P_i(c) + P_i(c)-a_i\|_{2,\mathscr{T}} \leq \|P_i(Q(a,b))-P_i(c)\|_{2,\mathscr{T}} + \|P_i(c)-a_i\|_{2,\mathscr{T}} \\ \leq L'(E \cap I_2) \cdot \|Q(a,b)-(c,d)\|_{2,E \cap I_2,\mathscr{T}} + \|P_i(c)-a_i\|_{2,\mathscr{T}} < L'(E \cap I_2) \cdot \eta \cdot |E \cap I_2|^{1/2} + \epsilon < 2\epsilon.
\end{align*}
Indeed, the $P_i$ consist only of variables indexed by $I_2$. A slight abuse of notation must be pointed out. The evaluation $Q(a,b)$ is an $I_2 \sqcup J_2$ tuple, while the Lipschitz constants $L'$ were defined on $S$-bounded $I_2$ tuples. However, since the $P_i$ consist only of variables indexed by $I_2$, in applying the $P_i$ one may view $Q(a,b)$ as an $S$-bounded $I_2$ tuple by considering only the $I_2$ indices in $Q(a,b)$. 

\medskip

Now, consider the following weak-$*$ neighborhood of $L_{(a,b)}$ defined by:
\begin{align*}
    U_0 \defeq \{\ell \in \Sigma_{R,I_1 \sqcup J_1} \, | \, \ell[(P_i \circ Q -T_i)^*(P_i \circ Q-T_i)]^{1/2} < 2\epsilon \, \, \text{for all} \, i \in F\}.
\end{align*}
As the intersection of finitely many opens indexed by $F$, the set $U_0$ is open. Moreover, the above estimate demonstrates that $U_0$ indeed contains $L_{(a,b)}$ as the $\|\cdot\|_{2,\mathscr{T}}$ norm is defined uniformly over all traces on $\mathscr{B}$. By definition, each microstate $(A,B)$ belonging to $\Gamma_R^{(N)}(U_0)$ satisfies:
\begin{align*}
    \|P_i(Q(A,B))-A_i\|_2 < 2\epsilon
\end{align*}
for all $i$ in $F$, where $\|\cdot\|_2$ is the $2$-norm on $N \times N$ matrices arising from the normalized trace. Applying (1) of the above Lemma, there exists a weak-$*$ neighborhood $U_1$ of $L_{(a,b)}$ such that:
\begin{align*}
    Q(\Gamma_R^{(N)}(U_1)) \subset \Gamma_S^{(N)}(\mathcal{V}).
\end{align*}
By design, $\Gamma_S^{(N)}(\mathcal{V}) \subset \Gamma_S^{(N)}(\mathcal{O})$, so that $Q(\Gamma_R^{(N)}(U_1)) \subset \Gamma_S^{(N)}(\mathcal{O})$. Consider now the intersection $U_0 \cap U_1$ and fix a microstate $(A,B)$ belonging to $\Gamma_R^{(N)}(U_0 \cap U_1)$. Note that $P_i$ is $\|\cdot\|_{2,E \cap I_2}$-$\|\cdot\|_2$ Lipschitz on $(S\text{Ball}(M_N(\C)))^{I_2 \sqcup J_2}$ for each $i$ in $F$, where:
\begin{align*}
    \|X\|_{2,E \cap I_2} \defeq \bigg( \sum_{j \in E \cap I_2} \|X_j\|_2^2\bigg)^{1/2}.
\end{align*}
Moreover, the Lipschitz constants are independent of $N$ (see the above remark). Denote the maximum Lipschitz constant among the $P_i$ by $L$. Let $\Omega_N$ be $(E \cap I_2, \epsilon/(|E \cap I_2|^{1/2}(L+1))$ orbit dense in $\Gamma_S^{(N)}(\mathcal{O})$ with:
\begin{align*}
    K_{E \cap I_2,\epsilon/(|E \cap I_2|^{1/2}(L+1))}^\text{orb}(\Gamma_S^{(N)}(\mathcal{O})) = |\Omega_N|. 
\end{align*}
Now, recall that $\|P_i(Q(A,B))-A_i\|_2 < 2\epsilon$ for each $i$ in $F$ as $(A,B)$ belongs to $\Gamma_R^{(N)}(U_0)$. By the definition of $\Omega_N$, there is a tuple $Y$ in $\Omega_N$ and a unitary $V$ for which:
\begin{align*}
    \|Q_j(A,B)-VY_jV^*\|_2 < \frac{\epsilon}{|E \cap I_2|^{1/2}(L+1)}. 
\end{align*}
for each $j$ belonging to $E \cap I_2$. Indeed, $(A,B)$ also belongs to $\Gamma_R^{(N)}(U_1)$. Therefore,
\begin{align*}
    \|Q(A,B)-VYV^*\|_{2,E \cap I_2} < |E \cap I_2|^{1/2} \cdot \frac{\epsilon}{|E \cap I_2|^{1/2}(L+1)} = \frac{\epsilon}{L+1}. 
\end{align*}
Since both $Q(A,B)$ and $Y$ are $S$-bounded tuples, one has for each member $i$ of $F$ the following:
\begin{align*}
    \|P_i(Q(A,B))-P_i(VYV^*)\|_2 \leq L \cdot  \|Q(A,B)-VYV^*\|_{2,E \cap I_2} < L \cdot \frac{\epsilon}{L+1} < \epsilon.
\end{align*}
Hence, appealing to the triangle inequality, the following holds for each index $i$ in $F$:
\begin{align*}
    \|A_i - P_i(VYV^*)\|_2 = \|A_i - P_i(Q(A,B))+P_i(Q(A,B))-P_i(VYV^*)\|_2 \\ \leq \|A_i-P_i(Q(A,B))\|_2 + \|P_i(Q(A,B))-P_i(VYV^*)\|_2 < 3\epsilon.
\end{align*}
That is, $\Gamma_R^{(N)}(U_0 \cap U_1)$ is $(F,3\epsilon)$ contained in the unitary conjugation orbit of $P(\Omega_N)$. Scaling $\epsilon$ by a constant $\kappa$ independent of $N$, there exists an $(F,\kappa \epsilon)$ orbit dense subset of $\Gamma_R^{(N)}(U_0 \cap U_1)$. The constant $\kappa$ arises as orbit dense collections are defined as subsets of the set to be packed. Therefore,
\begin{align*}
    \inf_{\mathscr{V} \supset L_{(a,b)}} \limsup_{N \rightarrow \infty} \frac{1}{N^2}\log K_{F,\kappa \epsilon}^\text{orb}(\Gamma_R^{(N)}(\mathscr{V})) \leq \limsup_{N \rightarrow \infty} \frac{1}{N^2} \log K_{E \cap I_2,\epsilon/[|E \cap I_2|^{1/2}(L+1)]}^\text{orb}(\Gamma_S^{(N)}(\mathcal{O})).
\end{align*}
Infimizing over $\mathcal{O}$, taking the supremum over finite subsets $F$ of $I_1$, and then taking the supremum over finite subsets $E'$ of $I_2$ (of which $E \cap I_2$ is an example), one obtains:
\begin{align*}
    h(a:b) \leq h(c:d).
\end{align*}
By symmetry, the reverse inequality holds and the proof is complete. 

\end{proof}

Now that invariance has been proved, one simply denotes the $1$-bounded entropy of $(M,\mathscr{T})$ by $h(M,\mathscr{T})$. Or, when the collection of traces is understood, simply $h(M)$. 

\section{The Variational Principle}

In this section, a result relating the $1$-bounded entropy developed above for $C^*$-algebras to the already existing $1$-bounded entropy for von Neumann algebras is proved. This will be the chief tool used in examples and applications below for the computation of the $1$-bounded entropy.

\medskip

First, consider a $C^*$-algebra $\mathscr{A}$. Now, set $\overline{\mathscr{A}}^{\mathscr{T}(\mathscr{A})}$ be the tracial completion of $\mathscr{A}$ with respect to the set $\mathscr{T}(\mathscr{A})$ of all traces on $\mathscr{A}$. For each such trace $\tau$, denote by $M_\tau$ the von Neumann algebra generated by the GNS representation $(\pi_\tau,L^2(\mathscr{A},\tau))$. 

\medskip

In the following Theorem, the case to keep in mind is the situation in which $K$ is the set $L_a$ of laws arising from all traces on a $C^*$ algebra $\mathscr{A}$ where $a$ is an $R$-bounded $I$ tuple which generates $\mathscr{A}$. 

\medskip

\noindent \textbf{Theorem 6.1 (Variational Principle):} Fix a compact, convex subset $K$ of the space $\Sigma_{R,I}$ of laws. One has the equality:
\begin{align*}
    h(K) = \sup_{\tau \in K} h(\{\ell_\tau\}) =\sup_{\tau \in K} h(M_\tau).
\end{align*}
 On the right hand side, $h$ refers to the $1$-bounded entropy of a tracial von Neumann algebra. It is defined by $h(M_\tau) \defeq h(\{\ell_\tau\})$. Henceforth, it will be convenient to write $h(\ell) \defeq h(\{\ell\})$.
\begin{proof}
    Consider a compact, convex subset $K$ of $\Sigma_{R,I}$. Since $\{\ell\} \subset K$ for each law $\ell$ in $K$, one has by definition:
    \begin{align*}
        h(\ell) \leq h(K).
    \end{align*}
    Therefore, $\sup_{\ell \in K} h(\ell) \leq h(K)$. Conversely, fix $\epsilon>0$ and a finite subset $E$ of the index set $I$. Fix a decreasing net $(\mathscr{V}_j)_{j \in J}$ of weak-$*$ open neighborhoods of the origin within the algebraic dual of $\C^*\langle T_i \, | \, i \in I \rangle$. Assume additionally that
    \begin{align*}
        \bigcap_{j \in J} \mathscr{V}_j = \{0\},
    \end{align*}
    and that $\mathscr{V}_j = - \mathscr{V}_j$ for each $j$. That is, the $\mathscr{V}_j$ are symmetric open neighborhoods of the origin which decrease to $0$. Fix another decreasing family $\mathscr{W}_j \subset \mathscr{V}_j$ of convex, symmetric, open neighborhoods of $0$ with the property that $10 \mathscr{W}_j \subset \mathscr{V}_j$ for each $j$. Fix $j$ and observe that:
    \begin{align*}
        K \subset \bigcup_{\ell \in K} [(\ell + \mathscr{W}_j) \cap \Sigma_{R,I}].
    \end{align*}
    Since $K$ is compact, there exists a finite subset $F_j$ of $K$ such that:
    \begin{align*}
        K \subset \bigcup_{\ell \in F_j} [(\ell + \mathscr{W}_j) \cap \Sigma_{R,I}].
    \end{align*}
    In order to obtain information about microstates, one needs to have an open neighborhood of $K$. For this, simply consider the sumset of $K$ and $\mathscr{W}_j$. One has:
    \begin{align*}
        (K+\mathscr{W}_j) \cap \Sigma_{R,I} \subset \bigcup_{\ell \in F_j} [(\ell + 2\mathscr{W}_j) \cap \Sigma_{R,I}].
    \end{align*}
    Indeed, $\mathscr{W}_j+\mathscr{W}_j \subset 2\mathscr{W}_j$ by convexity of $\mathscr{W}_j$. Now, this yields the following containment of the associated microstates spaces for each $N$:
    \begin{align*}
        \Gamma_R^{(N)}((K+\mathscr{W}_j) \cap \Sigma_{R,I}) \subset \bigcup_{\ell \in F_j} \Gamma_R^{(N)}((\ell + 2\mathscr{W}_j) \cap \Sigma_{R,I}).
    \end{align*}
    Taking orbital covering numbers, this gives the following inequality:
    \begin{align*}
        K_{E,\epsilon}^\text{orb}(\Gamma_R^{(N)}((K+\mathscr{W}_j) \cap \Sigma_{R,I})) \leq \sum_{\ell \in F_n} K_{E,\epsilon/2}^\text{orb}(\Gamma_R^{(N)}((\ell + 2\mathscr{W}_j) \cap \Sigma_{R,I})).
    \end{align*}
    Indeed, it is evident that the orbital covering number is subadditive under unions, where the $\epsilon/2$ adjustment comes from the fact that the centers in the packing must belong to the set to be packed. Now, observe:
    \begin{align*}
       \log K_{E,\epsilon}^\text{orb}(\Gamma_R^{(N)}((K+\mathscr{W}_j))) \leq \log \bigg(\sum_{\ell \in F_j} K_{E,\epsilon/2}^\text{orb}(\Gamma_R^{(N)}((\ell + 2\mathscr{W}_j) \cap \Sigma_{R,I}) )\bigg) \leq \log |F_j| \max_{\ell \in F_j}K_{E,\epsilon/2}^\text{orb}(\Gamma_R^{(N)}((\ell + 2\mathscr{W}_j) \cap \Sigma_{R,I})).
    \end{align*}
    Hence,
    \begin{align*}
        \limsup_{N \rightarrow \infty}\frac{1}{N^2} \log K_{E,\epsilon}^\text{orb}(\Gamma_R^{(N)}((K+\mathscr{W}_j))) \leq \limsup_{N \rightarrow \infty} \frac{1}{N^2}\log[ |F_j| \max_{\ell \in F_j}K_{E,\epsilon/2}^\text{orb}(\Gamma_R^{(N)}((\ell + 2\mathscr{W}_j) \cap \Sigma_{R,I}))] \\ \leq \max_{\ell \in F_j} \limsup_{N \rightarrow \infty} \frac{1}{N^2} \log K_{E,\epsilon/2}^\text{orb}(\Gamma_R^{(N)}((\ell + 2\mathscr{W}_j) \cap \Sigma_{R,I})).
    \end{align*}
    Denote by $\ell_j$ the law in $F_j$ for which the maximum on the right hand side is achieved. Moreover, denote by $\mathcal{O}_j$ the neighborhood $(K+\mathscr{W}_j) \cap \Sigma_{R,I}$. The above computations demonstrate the following:
    \begin{align*}
        \limsup_{N \rightarrow \infty} \frac{1}{N^2} \log K_{E,\epsilon}^\text{orb}(\Gamma_R^{(N)}(\mathcal{O}_j)) \leq \limsup_{N \rightarrow \infty} \frac{1}{N^2} \log K_{E,\epsilon/2}^\text{orb}(\Gamma_R^{(N)}((\ell_j+2\mathscr{W}_j) \cap \Sigma_{R,I})).
    \end{align*}
    
    Now, recall that $K$ is compact. Hence, perhaps passing to a subnet, without loss of generality suppose that $\ell_j \rightarrow \ell$ weak-$*$ for some law $\ell$ belonging to $K$. Now, since the $\mathscr{W}_j$ are open, for each index $j_0$ and for all $j$ sufficiently large one has that $\ell-\ell_j$ belongs to $\mathscr{W}_{j_0}$. Note (for fixed $j_0$ and for $j$ sufficiently large):
    \begin{align*}
      \ell_j+2\mathscr{W}_j \subset   \ell_j + 2\mathscr{W}_{j_0} \subset \ell +3\mathscr{W}_{j_0} \subset \ell + \mathscr{V}_{j_0}.
    \end{align*}
    Hence, 
    \begin{align*}
        \limsup_{N \rightarrow \infty} \frac{1}{N^2} \log K_{E,\epsilon}^\text{orb}(\Gamma_R^{(N)}(\mathcal{O}_j)) \leq \limsup_{N \rightarrow \infty} \frac{1}{N^2} \log K_{E,\epsilon/2}^\text{orb}(\Gamma_R^{(N)}((\ell + \mathscr{V}_{j_0}) \cap \Sigma_{R,I})).
    \end{align*}
    Infimizing over $j_0$, one obtains:
    \begin{align*}
        \inf_{\mathcal{O} \supset K}  \limsup_{N \rightarrow \infty} \frac{1}{N^2} \log K_{E,\epsilon}^\text{orb}(\Gamma_R^{(N)}(\mathcal{O})) \leq \inf_{\mathcal{S} \supset \{\ell\}}\limsup_{N \rightarrow \infty} \frac{1}{N^2} \log K_{E,\epsilon/2}^\text{orb}(\Gamma_R^{(N)}(\mathcal{S})).
    \end{align*}
    Indeed, exchanging the infimum over $j_0$ with the infimum over all weak-$*$ neighborhoods of $K$ only decreases the left hand side further. On the other hand, since the $\mathscr{V}_{j_0}$ shrink to $0$, replacing the infimum over $j_0$ on the right hand side with the infimum over all neighborhoods of $\ell$ does not change the right hand side. Now, one may not simply take the supremum over $\epsilon$ and $E$, because $\ell$ implicitly depends on these by construction. Therefore, taking first the supremum over all laws yields:
    \begin{align*}
        \inf_{\mathcal{O} \supset K}  \limsup_{N \rightarrow \infty} \frac{1}{N^2} \log K_{E,\epsilon}^\text{orb}(\Gamma_R^{(N)}(\mathcal{O})) \leq \sup_{\ell \in K}\inf_{\mathcal{S} \supset \{\ell\}}\limsup_{N \rightarrow \infty} \frac{1}{N^2} \log K_{E,\epsilon/2}^\text{orb}(\Gamma_R^{(N)}(\mathcal{S})).
    \end{align*}
    Now that the dependence on $E$ and $\epsilon$ washed out in the supremum, one can take the supremum over $E$ and $\epsilon$. Taking the supremum over $E$ and $\epsilon$ and swapping it with the supremum over $K$ on the right hand side yields:
    \begin{align*}
        h(K) \leq \sup_{\ell \in K} h(\ell)
    \end{align*}
    as desired.
\end{proof}

\noindent \textbf{Remark:} In work of Jekel \cite{JekelVar}, a similar variational principle arises relating the entropy of an existential or quantifier free type to the entropies of each of the full types compatible with it. 

\section{Topological $1$-Bounded Entropy}

In the previous section, a notion of $1$-bounded entropy for a $C^*$ algebra $\mathscr{A}$ was introduced in which one takes the $1$-bounded entropy of the tracial completion. In this section, a distinct notion of the entropy involving the operator norm is introduced. This method takes into account recent developments in strong convergence. It will later be demonstrated that this topological $1$-bounded entropy is bounded above by the $1$-bounded entropy of the previous section.

\medskip

Fix a $C^*$-algebra $\mathscr{A}$ and let $x = (x_i)_{i \in I}$ be a generating tuple for $\mathscr{A}$. Assume that $C>0$ is such that $\|x_i\| \leq C$ for all $x$. Fix $\delta >0$. Now, for each finite subset $Q$ of $\C^*\langle T_i \, | \, i \in I \rangle$ consider the following topological microstates space:
\begin{align*}
    \Gamma_C^{(N,\text{top})}(x;Q,\delta) \defeq \bigcap_{P \in Q} \{A \in (C\text{Ball}M_N(\C))^I \, | \, \|P(A)\| \leq \|P(x)\| + \delta\}. 
\end{align*}
Here, $\|\cdot\|$ refers to both the operator norm on matrices and the $C^*$-norm on $\mathscr{A}$. It will be evident from context which is being used. Now, consider a subset $\Omega$ of $M_N(\C)^I$. For each tuple $B$ in $\Omega$, each finite subset $F \subset I$, and each $\epsilon >0$, one defines the $(F,\epsilon)$ orbital neighborhood $N_{F,\epsilon}^\text{orb}(B,\|\cdot\|)$ with respect to the operator norm as follows:
\begin{align*}
    N_{F,\epsilon}^\text{orb}(B,\|\cdot\|) \defeq \{C \in M_N(\C)^I \, | \, \exists U \in U(N) \, \text{such that} \, \|C_i-UB_iU^*\| < \epsilon \, \forall i \in F\}. 
\end{align*}
Define further
\begin{align*}
    N_{F,\epsilon}^\text{orb}(\Omega,\|\cdot\|) = \bigcup_{B \in \Omega}N_{F,\epsilon}^\text{orb}(B,\|\cdot\|).
\end{align*}
A subset $\Omega_0$ of $\Omega$ is called $(F,\epsilon)$ orbit dense with respect to the operator norm provided:
\begin{align*}
    \Omega \subset N_{F,\epsilon}(\Omega_0,\|\cdot\|).
\end{align*}
One defines the $(F,\epsilon)$ orbital covering number with respect to the operator norm $K_{F,\epsilon}^\text{orb}(\Omega,\|\cdot\|)$ to be the minimal cardinality of an $(F,\epsilon)$ orbit dense subset $\Omega_0$ of $\Omega$. By compactness, this quantity is finite provided $\Omega$ is contained in a product of closed operator norm balls. This is the case for the topological microstates spaces due to the norm cutoff $C>0$. This yields the following definition. As in the above case, it is convenient later on to work with an "in the presence" definition of the entropy.

\medskip

\noindent \textbf{Definition 7.1 (Topological $1$-Bounded Entropy):} Fix $C>0$ and let $\mathscr{A} \subset \mathscr{B}$ be an inclusion of $C^*$ algebras. Suppose that $\mathscr{A}$ is generated by a $C$-bounded tuple $a=(a_i)_{i \in I}$, and suppose additionally that $b=(b_j)_{j \in J}$ is a $C$-bounded tuple such that $a$ and $b$ jointly generate $\mathscr{B}$. Define
\begin{align*}
    h^\text{top}(a:b) = \sup_{\epsilon, F} \inf_{\delta,Q} \limsup_{N \rightarrow \infty} \frac{1}{N^2}\log K_{F,\epsilon}^\text{orb}(\Gamma_C^{(N,\text{top})}((a,b);Q,\delta),\|\cdot\|),
\end{align*}
where $\epsilon>0$, the set $F$ is a finite collection of indices in $I$, the set $Q$ is a finite collection of $*$-polynomials in $\C^*\langle T_i \, | \, i \in I \sqcup J \rangle$, and $\delta >0$. This quantity is referred to as the topological $1$-bounded entropy of $\mathscr{A}$ in the presence of $\mathscr{B}$.

\medskip

Now, a proof of invariance for the topological $1$-bounded entropy is presented. After this is carried out, the quantity may simply be denoted $h^\text{top}(\mathscr{A}:\mathscr{B})$. First, a Lemma is required. Its proof is omitted, and amounts to an induction on the degree and the triangle inequality. Indeed, it suffices to verify the claim on monomials. 

\medskip

\noindent \textbf{Lemma 7.1:} Fix $R>0$ and $I$ an index set. Consider a $*$-polynomial $P$ belonging to $\C^*\langle T_i \, | \, i \in I \rangle$. If $E$ is a finite subset of the index set $I$ which contains all indices $i$ corresponding to variables appearing nontrivially in a monomial of $P$, then there exists a constant $L>0$ such that
\begin{align*}
    \|P(x)-P(y)\| \leq L \cdot \max_{i \in E} \|x_i-y_i\|
\end{align*}
for all $R$-bounded tuples $x$ and $y$. In the sequel, the norm $\max_{i \in E} \|x_i\|$ will simply be denoted $\|x\|_E$. 

\medskip

\noindent \textbf{Theorem 7.1 (Invariance of Topological $1$-Bounded Entropy):} Let $\mathscr{A} \subset \mathscr{B}$ be an inclusion of $C^*$ algebras. Fix norm cutoffs $R,S>0$ and index sets $I_1,J_1,I_2,J_2$. Suppose $a=(a_i)_{i \in I_1}$ and $b=(b_j)_{j \in J_1}$ are both $R$-bounded such that $a$ generates $\mathscr{A}$ and $a$ and $b$ jointly generate $\mathscr{B}$. Similarly, suppose $c=(c_i)_{i \in I_2}$ and $d=(d_j)_{j \in J_2}$ are $S$-bounded such that $c$ generates $\mathscr{A}$ and $c$ and $d$ jointly generate $\mathscr{B}$. Then
\begin{align*}
    h^\text{top}(a:b) = h^\text{top}(c:d).
\end{align*}
\begin{proof}
    Fix a finite subset $F$ of $I_1$ and $\epsilon>0$. Since $c$ generates $\mathscr{A}$, there exist polynomials $(P_i)_{i \in F}$ belonging to $\C^* \langle T_j \, | \, j \in I_2 \rangle$ such that:
    \begin{align*}
        \|P_i(c)-a_i\| < \epsilon
    \end{align*}
    for each index $i$ in $F$. Additionally, one may take $\|P_i\|_S \leq R$ for each $i$ in $F$. For $i$ a member of $I_1 \sqcup J_1 \backslash F$, define $P_i \equiv 0$ such that $P_i(c)$ may be considered an $R$-bounded $I_1 \sqcup J_1$ tuple. 

    \medskip

    Now, fix a finite collection $Q$ of $*$-polynomials belonging to $\C^* \langle T_j \, | \, j \in I_2 \sqcup J_2 \rangle$ and $\delta>0$. Let $E$ be a finite subset of $I_2 \sqcup J_2$ which contains all indices of variables $T_j$ or $T_j^*$ appearing nontrivially in members of $Q$ and within the $P_i$ for $i$ an index in $F$. Denote the largest $\|\cdot\|_E$-$\|\cdot\|$ Lipschitz constant among the members of $Q$ and the $P_i$ on $S$-bounded tuples by $L$. Such a constant exists due to the above Lemma.

    \medskip

    Fix $0<\hat{\delta} \ll \delta$ with $\hat{\delta}<\epsilon$, and choose $\eta>0$ such that $L \eta < \epsilon$ and $L \eta < \delta-\hat{\delta}$. Since $(a,b)$ generates $\mathscr{B}$, there exist polynomials $Z_j$ in $\C^*\langle T_i \, | \, i \in I_1 \sqcup J_1 \rangle$ with $\|Z_j\|_R \leq S$ and:
    \begin{align*}
        \|Z_j(a,b)-(c,d)_j\| < \eta
    \end{align*}
    for each index $j$ belonging to $E$. For $j$ a member of $I_2 \sqcup J_2 \backslash E$, define $Z_j \equiv 0$, so that $Z(a,b)$ may be considered an $S$-bounded $I_2 \sqcup J_2$ tuple. Now, set
    \begin{align*}
        \hat{Q} \defeq \{X \circ Z_j \, | \, X \in Q\} \cup \{P_i \circ Z - T_i \, | \, i \in F\}
    \end{align*}
    and consider the associated topological microstates space $\Gamma_R^{(N,\text{top})}((a,b);\hat{Q},\hat{\delta})$. Observe that for each microstate $(A,B)$ belonging to this space and $i$ a member of $F$ there holds:
    \begin{align*}
        \|P_i(Z(A,B))-A_i\| \leq \|P_i(Z(a,b))-a_i\| + \hat{\delta} \leq \|P_i(Z(a,b))-P_i(c)\| + \|P_i(c)-a_i\| \\ \leq L \cdot \|Z(a,b)-(c,d)\|_E + 2\epsilon < L\eta + 2\epsilon < 3\epsilon.
    \end{align*}
    Indeed, $P_i(c,d)=P_i(c)$ as the $P_i$ involve only indices from $I_2$, one has $\hat{\delta} < \epsilon$, and the $P_i(c)$ approximate the $a_i$ within $\epsilon$. Now, it will be demonstrated that $Z(\Gamma_R^{(N,\text{top})}((a,b);\hat{Q},\hat{\delta}))$ is contained within $\Gamma_S^{(N,\text{top})}((c,d);Q,\delta)$. Fix a microstate $(A,B)$ belonging to $\Gamma_R^{(N,\text{top})}((a,b);\hat{Q},\hat{\delta})$ and observe the following for each $X$ in $Q$ by the reverse triangle inequality:
    \begin{align*}
        \|X(Z(A,B))\| \leq \|X(Z(a,b))\| + \hat{\delta} \leq | \,  \|X(Z(a,b))\|-\|X(c,d)\| \,  | + \|X(c,d)\| + \hat{\delta} \\  \leq L \cdot \|Z(a,b)-(c,d)\|_E + \|X(c,d)\| + \hat{\delta} < L\eta + \|X(c,d)\| + \hat{\delta} < \|X(c,d)\| + \delta.
    \end{align*}
    Indeed, the $X \circ Z$ belong to $\hat{Q}$ and $L \eta < \delta - \hat{\delta}$. By definition, the above estimate demonstrates that $Z(A,B)$ belongs to $\Gamma_S^{(N,\text{top})}((b,c);Q,\delta)$. 

    \medskip
    
    Now, note that each $P_i$ for $i$ an index in $F$ is $\|\cdot\|_{E \cap I_2}$-$\|\cdot\|$ Lipschitz on $(S\text{Ball}(M_N(\C)))^{I_2}$ with Lipschitz constant independent of $N$. Denote the maximum Lipschitz constant among the $P_i$ by $L'$. Let $\Omega_N \subset \Gamma_S^{(N,\text{top})}((b,c);Q,\delta)$ be $(E \cap I_2,\epsilon/(L'+1))$ orbit dense with:
    \begin{align*}
        K_{E \cap I_2,\epsilon/(L'+1)}^\text{orb}(\Gamma_S^{(N,\text{top})}((b,c);Q,\delta),\|\cdot\|) = |\Omega_N|.
    \end{align*}
    Fix a microstate $(A,B)$ belonging to $\Gamma_R^{(N,\text{top})}((a,b);\hat{Q},\hat{\delta})$ and recall that
    \begin{align*}
        \|P_i(Z(A,B))-A_i\| <3\epsilon
    \end{align*}
    for each $i$ in $F$. Moreover, since $Z(A,B)$ is a member of $\Gamma_S^{(N,\text{top})}((b,c);Q,\delta)$, there exists a tuple $Y$ belonging to $\Omega_N$ and a unitary $V$ such that:
    \begin{align*}
        \|Z_j(A,B)-VY_jV^*\| < \frac{\epsilon}{L'+1}
    \end{align*}
    for each $j$ in $E \cap I_2$. Therefore,
    \begin{align*}
        \|Z(A,B)-VYV^*\|_{E \cap I_2} < \frac{\epsilon}{L'+1}.
    \end{align*}
    Since both $Z(A,B)$ and $Y$ are $S$-bounded, one has for each index $i$ in $F$ the following:
    \begin{align*}
        \|P_i(Z(A,B))-P_i(VYV^*)\| \leq L' \cdot \|Z(A,B)-VYV^*\|_{E \cap I_2} < L' \cdot \frac{\epsilon}{L'+1} < \epsilon.
    \end{align*}
    This in turn yields the following for each index $i$ in $F$:
    \begin{align*}
        \|A_i-P_i(VYV^*)\| \leq \|A_i-P_i(Z(A,B))\| + \|P_i(Z(A,B))-P_i(VYV^*)\| < 3\epsilon + \epsilon = 4\epsilon.
    \end{align*}
    That is, $\Gamma_R^{(N,\text{top})}((a,b);\hat{Q},\hat{\delta})$ is $(F,4\epsilon)$ contained in the unitary conjugation orbit of $P(\Omega_N)$ with respect to the operator norm. Indeed, unitary conjugation commutes with the application of $*$-polynomials. Scaling $\epsilon$ by $\kappa>0$ independent of $N$, there exists an $(F,\kappa \epsilon)$ orbit dense subset of $\Gamma_R^{(N,\text{top})}((a,b);\hat{Q},\hat{\delta})$ with respect to the operator norm of size $|\Omega_N|$. Hence,
    \begin{align*}
        \inf_{\Tilde{Q},\Tilde{\delta}} \limsup_{N \rightarrow \infty} \frac{1}{N^2}\log K_{F,\kappa \epsilon}^\text{orb}(\Gamma_R^{(N,\text{top})}((a,b);\Tilde{Q},\Tilde{\delta}),\|\cdot\|) \leq \limsup_{N \rightarrow \infty} \frac{1}{N^2} \log K_{E \cap I_2,\epsilon/(L'+1)}^\text{orb}(\Gamma_S^{(N,\text{top})}((b,c);Q,\delta),\|\cdot\|)
    \end{align*}
    Infimizing over $Q$ and $\delta$, taking the supremum over $E' \subset I_2$ (of which $E \cap I_2$ is an example), taking the supremum over $F \subset I_1$, and finally taking the supremum over $\epsilon>0$ one obtains:
    \begin{align*}
        h^\text{top}(a:b) \leq h^\text{top}(c:d).
    \end{align*}
    By symmetry, the reverse inequality holds and the claim is proved.
\end{proof}

Given that the topological $1$-bounded entropy is an invariant, it may be written $h^\text{top}(\mathscr{A}:\mathscr{B})$ for $\mathscr{A} \subset \mathscr{B}$ an inclusion of $C^*$-algebras. 

\section{Comparing The Different Notions}

In this section, it will be demonstrated that if $a$ is a tuple generating a $C^*$-algebra $\mathscr{A}$, then $h^\text{top}(a) \leq h(L_a)$. That is, the topological $1$-bounded entropy of $\mathscr{A}$ does not exceed the $1$-bounded entropy of the universal tracial completion of $\mathscr{A}$. A similar inequality holds between the "in the presence" versions of both entropies. The proof of the comparison proceeds by demonstrating a containment of microstates spaces, irrespective of how they are packed. Therefore, the "in the presence" version introduces no additional technical difficulties. Hence, for ease of exposition, the proofs in this section involve only the topological $1$-bounded entropy and the $1$-bounded entropy of a single $C^*$ algebra $\mathscr{A}$. 

\medskip

First, a Lemma.

\medskip

\noindent \textbf{Lemma 8.1:} Suppose that $\mathscr{A}$ is a $C^*$-algebra and fix an $R$-bounded generating tuple $a=(a_i)_{i \in I}$ of $\mathscr{A}$ where $R$ strictly exceeds the norm of each $a_i$. Suppose additionally that $\mathscr{A}$ is separable, and that $I$ is countable. Set $\mathcal{O} \subset \Sigma_{R,I}$ to be a weak-$*$ neighborhood of $L_a$. Then there exists $\delta>0$ and a finite set $E \subset \C^*\langle T_i \, | \, i \in  I \rangle$ of star polynomials such that if $(B,\tau)$ is a $C^*$-algebra equipped with a tracial state $\tau$ and $b$ is an $R$-bounded $I$-tuple in $B$ such that:
\begin{align*}
    \|Q(b)\| \leq \|Q(a)\| + \delta
\end{align*}
for each $Q$ belonging to $E$, then $\ell_b$ belongs to $\mathcal{O}$. Here, the law of $b$ is the composition of the evaluation map on $b$ and the tracial state $\tau$. 
\begin{proof}
    Suppose towards a contradiction that for all $\delta >0$ and for all finite subsets $E$ of $\C^*\langle T_i \, | \, i \in  I \rangle$ there exists a $C^*$-algebra $(B,\tau)$ (depending on $\delta$ and $E$) equipped with a tracial state $\tau$ and an $R$-bounded $I$-tuple $b$ of elements of $B$ satisfying
    \begin{align*}
    \|Q(b)\| \leq \|Q(a)\| + \delta
    \end{align*}
    for all members $Q$ of $E$, but such that $\ell_b$ does not belong to $\mathcal{O}$.

    \medskip
    
    Fix a dense subfield $K$ of $\C$, say $K=\Q(i)$. Enumerate the members of $K^*\langle T_i \, | \, i \in I \rangle$ by $\{P_k\}_{k=1}^\infty$. For each $n$, set 
    \begin{align*}
        E_n \defeq \{P_1,P_2,\dots,P_n\}. 
    \end{align*}
    Fix positive real numbers $\delta_n \rightarrow 0$. By assumption, there exists for each $n$ a $C^*$-algebra $(B_n,\tau_n)$ equipped with a tracial state $\tau_n$  and an $R$-bounded $I$-tuple $b_n$ of elements of $B_n$ satisfying:
    \begin{align*}
        \|Q(b_n)\| \leq \|Q(a)\| + \delta_n
    \end{align*}
    for each member $Q$ of $E_n$, and such that $\ell_{b_n}$ does not belong to $\mathcal{O}$. Observe that if $x$ is any $R$-bounded $I$-tuple in a $C^*$-algebra, then:
    \begin{align*}
        \bigcup_{n \geq 1} \{P(x) \, | \, P \in E_n\}
    \end{align*}
    is dense in $C^*(x)$. This coupled with the fact that $\delta_n \rightarrow 0$ and the observation that
    \begin{align*}
        \|Q(b_n)\| \leq \|Q(a)\| + \delta_n
    \end{align*}
    for all $Q$ in $E_n$ yields an embedding 
    \begin{align*}
       \iota: \mathscr{A} \hookrightarrow \prod_{n \rightarrow \omega} B_n
    \end{align*}
    of $\mathscr{A}$ within the operator norm ultraproduct of the $B_n$. Here, $\omega$ denotes a free ultrafilter on $\N$. Denote the image $\iota(a)$ by $b$. Note that for each star polynomial $P$, one has:
    \begin{align*}
        \ell_b(P) = \lim_{n \rightarrow \omega} \tau_n(P(b_n)) = \lim_{n \rightarrow \omega} \ell_{b_n}(P).
    \end{align*}
    That is, $\ell_{b_n} \rightarrow \ell_b$ in the weak-$*$ topology. This is by construction, as the ultraproduct is equipped with the trace $\lim_{n \rightarrow \omega} \tau_n$. Since $\ell_{b_n}$ does not belong to $\mathcal{O}$ for each $n$, neither does $\ell_b$. On the other hand, pulling back the trace $\lim_{n \rightarrow \omega} \tau_n$ under $\iota$ yields a trace $\tau$ on $\mathscr{A}$ so that
    \begin{align*}
        \ell_b(P) = \tau(P(a))
    \end{align*}
    for all $*$-polynomials $P$. Hence, $\ell_b$ belongs to $L_a$, a contradiction. 
\end{proof}

\medskip

The following Lemma is similar to Lemma 2.5 in Hayes \cite{Reg}. The objective is to demonstrate that the operator norm may also be used to compute the entropy. 

\medskip

\noindent \textbf{Lemma 8.2:} Fix $R>0$ and a natural number $N$ and suppose that $\Omega$ is a subset of $(R\text{Ball}(M_N(\C)))^I$ for some index set $I$. Fix a finite subset $F$ of $I$ and $\epsilon > 0$. Then there is a constant $C>0$ (depending only on $\epsilon$) such that for all $0<\delta<\epsilon$ one has:
\begin{align*}
   K_{F,2(2R+2)\epsilon}^\text{orb}(\Omega,\|\cdot\|) \leq K_{F,\delta}^\text{orb}(\Omega) \cdot \bigg[(1+N \frac{\delta^2}{\epsilon^2})^{} \bigg( \frac{C}{\epsilon} \bigg)^{2N^2(\delta^2/\epsilon^2)} \cdot \bigg(\frac{6R+\epsilon}{\epsilon}\bigg)^{2N^2(\delta^2/\epsilon^2)}\bigg]^{|F|}.
\end{align*}
\begin{proof}
    Set $D \subset \Omega$ to be $(F,\delta)$ orbit dense in $\Omega$ with respect to the $2$-norm. Let $\Xi$ be an $\epsilon$-dense subset with respect to the operator norm in the space $\mathscr{P}_{(\delta/\epsilon)^2}$ of projections in $M_N(\C)$ of trace not exceeding $(\delta/\epsilon)^2$. That is, for each $P$ in $\mathscr{P}_{(\delta/\epsilon)^2}$, there exists a $Q$ in $\Xi$ such that $\|P-Q\|<\epsilon$.

    \medskip
    
    Finally, for each projection $P$ in $\mathscr{P}_{(\delta/\epsilon)^2}$, set $D_P$ to be a maximal subset of $2R(P\text{Ball}(M_N(\C)),\|\cdot\|)$ such that $\|A-B\| \geq \epsilon$ for all distinct members $A$ and $B$ of $D_P$. By maximality, observe that $D_P$ is $\epsilon$-dense with respect to the operator norm in $2R(P\text{Ball}(M_N(\C)),\|\cdot\|)$. It will first be demonstrated that:
    \begin{align*}
      X \defeq  D + \bigg[ \bigcup_{Q \in \Xi} D_Q \bigg]^F
    \end{align*}
    contains $\Omega$ in its $(F,(2R+2)\epsilon)$ orbital neighborhood with respect to the operator norm. Note that a member of $X$ need not lie in $\Omega$. However, it is a straightforward application of the triangle inequality to see that one may produce from $X$ an honest $(F,2(2R+2)\epsilon)$ orbit dense \textit{subset} of $\Omega$ with the same cardinality as $X$. 

    \medskip

    To see that $X$ contains $\Omega$ in its $(F,(2R+2)\epsilon)$ orbital neighborhood with respect to the operator norm, fix a tuple $(A_i)_{i \in F}$ within $\Omega$. By definition, there exists a tuple $(Y_i)_{i \in F}$ belonging to $D$ and a unitary $U$ such that $\|UA_iU^*-Y_i\|_2 < \delta$ for each $i$ in $F$. For each $i$ in $F$, denote by $P_i$ the projection $\chi_{[\epsilon,\infty)}(|UA_iU^*-Y_i|)$. By construction:
    \begin{align*}
        \|(1-P_i)(UA_iU^*-Y_i)\| < \epsilon.
    \end{align*}
    Moreover, applying Markov's inequality:
    \begin{align*}
        \text{tr}(P_i) \leq \frac{1}{\epsilon^2}\text{tr}(|U A_i U^* - Y_i|^2) < \frac{\delta^2}{\epsilon^2}.
    \end{align*}
    Hence, for each $i$ in $F$ there is $Q_i$ in $\Xi$ with $\|P_i-Q_i\| < \epsilon$. Finally, for each $i$ in $F$ fix $C_i$ in $D_{Q_i}$ such that:
    \begin{align*}
        \|Q_i(UA_iU^*-Y_i)-C_i\| < \epsilon.
    \end{align*}
    By design, $(Y_i + C_i)_{i \in F}$ is a member of $X$. Observe the following for each $i$:
    \begin{align*}
        \|UA_i U^* - (Y_i+C_i)\| \leq \|(1-P_i)(UA_iU^*-Y_i)\| + \|P_i(UA_iU^*-Y_i)-C_i\| < \epsilon + \|P_i(UA_iU^*-Y_i)-C_i\| \\ \leq \epsilon + \|(P_i-Q_i)(UA_iU^*-Y_i)\| + \|Q_i(UA_iU^*-Y_i)-C_i\| \\ < \epsilon + \|P_i-Q_i\| \cdot \|UA_iU^*-Y_i\| + \epsilon \leq (2+2R)\epsilon.
    \end{align*}
   Thus $X$ contains $\Omega$ in its $(F,(2R+2)\epsilon)$ orbital neighborhood with respect to the operator norm. Therefore, the above remarks yield:
   \begin{align*}
       K_{F,2(2R+2)\epsilon}^\text{orb}(\Omega,\|\cdot\|) \leq |X|. 
   \end{align*}
   It remains to estimate the size of $X$. By the definition of $X$, one has:
   \begin{align*}
       |X| \leq |D| \cdot \bigg|\bigcup_{Q \in \Xi} D_Q \bigg|^{|F|} \leq |D| \cdot \bigg[ \sum_{Q \in \Xi} |D_Q| \bigg]^{|F|}. 
   \end{align*}
   Note that $|D|=K_{F,\delta}^\text{orb}(\Omega)$. Fix a projection $Q$ in $\Xi$. To control the size of $D_Q$, observe that the triangle inequality yields the following:
   \begin{align*}
       \bigg(2R+\frac{\epsilon}{3} \bigg) \text{Ball}(QM_N(\C),\|\cdot\|) \supset \bigcup_{C \in D_Q} C + \frac{\epsilon}{3}\text{Ball}(QM_N(\C),\|\cdot\|).
   \end{align*}
   Since $D_Q$ is $\epsilon$ separated by definition, the union on the right hand side is disjoint. Taking the volume on both sides:
   \begin{align*}
       \bigg(2R+\frac{\epsilon}{3}\bigg)^{2N^2\text{tr}(Q)} \geq |D_Q| \cdot \bigg(\frac{\epsilon}{3}\bigg)^{2N^2\text{tr}(Q)}. 
   \end{align*}
   Since the trace of $Q$ does not exceed $(\delta/\epsilon)^2$, the following holds:
   \begin{align*}
       |D_Q| \leq \bigg(\frac{6R+\epsilon}{\epsilon}\bigg)^{2N^2(\delta^2/\epsilon^2)}. 
   \end{align*}
   Substituting this into the above, one obtains:
   \begin{align*}
       |X| \leq K_{F,\delta}^\text{orb}(\Omega) \cdot \bigg[ |\Xi| \cdot \bigg(\frac{6R+\epsilon}{\epsilon}\bigg)^{2N^2(\delta^2/\epsilon^2)}  \bigg]^{|F|} = K_{F,\delta}^\text{orb}(\Omega) \cdot |\Xi|^{|F|} \cdot \bigg(\frac{6R+\epsilon}{\epsilon}\bigg)^{2N^2(\delta^2/\epsilon^2)|F|}.
   \end{align*}
   Now, applying the packing estimate of Szarek (see Lemma 3.2 of Hayes, Jekel, and Kunnawalkam Elayavalli \cite{PropT} for the precise formulation used here) one obtains a constant $C>0$ such that:
   \begin{align*}
       |\Xi| \leq (1+N \frac{\delta^2}{\epsilon^2}) \bigg( \frac{C}{\epsilon} \bigg)^{2N^2(\delta^2/\epsilon^2)}. 
   \end{align*}
   Putting it all together:
   \begin{align*}
        K_{F,2(2R+2)\epsilon}^\text{orb}(\Omega,\|\cdot\|) \leq |X| \leq K_{F,\delta}^\text{orb}(\Omega) \cdot (1+N \frac{\delta^2}{\epsilon^2})^{|F|} \bigg( \frac{C}{\epsilon} \bigg)^{2N^2(\delta^2/\epsilon^2)|F|} \cdot \bigg(\frac{6R+\epsilon}{\epsilon}\bigg)^{2N^2(\delta^2/\epsilon^2)|F|}.
   \end{align*}
   This completes the proof.
\end{proof}

With this Lemma in hand, one can show that $h$ may be computed via an operator norm packing. First, for a compact, convex subset $K$ of $\Sigma_{R,I}$, denote:
\begin{align*}
    h(K,\|\cdot\|) \defeq \sup_{F,\epsilon} \inf_{\mathcal{O} \supset K} \limsup_{N \rightarrow \infty} \frac{1}{N^2}\log (K_{F,\epsilon}^\text{orb}(\Gamma_R^{(N)}(\mathcal{O}),\|\cdot\|)).
\end{align*}
The formula is identical to that for $h(K)$, only the operator norm is used to compute the packing numbers of microstates spaces. The following is analogous to Corollary 2.6 of Hayes \cite{Reg}.

\medskip

\noindent \textbf{Proposition 8.1:} Suppose that $K$ is a compact, convex subset of $\Sigma_{R,I}$. Then
\begin{align*}
    h(K,\|\cdot\|) = h(K).
\end{align*}
\begin{proof}
    Evidently, $h(K) \leq h(K,\|\cdot\|)$, as $2$-norm balls contain operator norm balls of the same radius (the normalized trace is used to define the $2$-norm). Conversely, fix a weak-$*$ neighborhood $\mathcal{O}$ of $K$, a finite set $F$ of $I$, and $\epsilon>0$ . Observe the following from the Lemma for $0<\delta<\epsilon$:
    \begin{align*}
        \frac{1}{N^2} \log K_{F,2(2R+2)\epsilon}^\text{orb}(\Gamma_R^{(N)}(\mathcal{O}),\|\cdot\|) \leq \frac{1}{N^2} \log \bigg[ K_{F,\delta}^\text{orb}(\Gamma_R^{(N)}(\mathcal{O})) \cdot (1+N \frac{\delta^2}{\epsilon^2}) \bigg( \frac{C}{\epsilon}\bigg)^{2N^2(\delta^2/\epsilon^2)} \cdot \bigg( \frac{6R+\epsilon}{\epsilon} \bigg)^{2N^2(\delta^2/\epsilon^2)} \bigg]^{|F|}. 
    \end{align*}
    Taking the limit superior as $N \rightarrow \infty$ and infimizing over $\mathcal{O}$ yields the following:
    \begin{align*}
        \inf_{\mathcal{O} \supset K} \limsup_{N \rightarrow \infty}\frac{1}{N^2} \log K_{F,2(2R+2)\epsilon}^\text{orb}(\Gamma_R^{(N)}(\mathcal{O}),\|\cdot\|) \leq \inf_{\mathcal{O} \supset K} \limsup_{N \rightarrow \infty}\frac{1}{N^2} \log K_{F,\delta}^\text{orb}(\Gamma_R^{(N)}(\mathcal{O})) + 2|F|(\delta^2/\epsilon^2)[\log(C/\epsilon)+\log \bigg( \frac{6R+\epsilon}{\epsilon} \bigg)].
    \end{align*}
    Taking the limit superior as $\delta \rightarrow 0$, and then the supremum over $\epsilon > 0$, and finally taking the supremum over finite subsets $F$ of $I$ yields the claim. 
\end{proof}

\medskip

Now, the main result comparing the two invariants may be established.

\medskip

\noindent \textbf{Theorem 8.1:} Let $\mathscr{A}$ be a $C^*$-algebra and fix an $R$-bounded generating tuple $(a_i)_{i \in I}$, where $R$ is strictly larger than the norm of each $a_i$. Suppose additionally that $\mathscr{A}$ is separable, and that $I$ is countable. Then:
\begin{align*}
  h^\text{top}(\mathscr{A}) =  h^\text{top}(a) \leq h(L_a) = h(\mathscr{A}).
\end{align*}
\begin{proof}
    By the above proposition, it suffices to demonstrate that $h^\text{top}(a) \leq h(L_a,\|\cdot\|)$. Fix a weak-$*$ neighborhood $\mathcal{O}$ of $L_a$ within $\Sigma_{R,I}$. By the Lemma at the beginning of this section, there exists $\delta>0$ and a finite subset $E$ of $\C^*\langle T_i \, | \, i \in I \rangle$ such that:
    \begin{align*}
        \Gamma_R^{(N,\text{top})}(a;E,\delta) \subset \Gamma_R^{(N)}(\mathcal{O})
    \end{align*}
    for all $N$. Therefore, for each finite subset $F$ of the index set $I$, there holds:
    \begin{align*}
    \limsup_{N \rightarrow \infty} \frac{1}{N^2} \log K_{F,\epsilon}^\text{orb}( \Gamma_R^{(N,\text{top})}(a;E,\delta),\|\cdot\|) \leq  \limsup_{N \rightarrow \infty} \frac{1}{N^2} \log K_{F,\epsilon/2}^\text{orb}(\Gamma_R^{(N)}(\mathcal{O}),\|\cdot\|).
    \end{align*}
    The $\epsilon/2$ arises as the relevant correction for the fact that the covering numbers defined above are not monotone with respect to inclusions of sets. Indeed, the covering balls are required to have centers within the set to be covered. Infimizing over $\delta>0$ and finite subsets $E$ of $\C^*\langle T_i \, | \, i \in I \rangle$ yields:
    \begin{align*}
       \inf_{E,\delta} \limsup_{N \rightarrow \infty} \frac{1}{N^2} \log K_{F,\epsilon}^\text{orb}( \Gamma_R^{(N,\text{top})}(a;E,\delta),\|\cdot\|) \leq  \limsup_{N \rightarrow \infty} \frac{1}{N^2} \log K_{F,\epsilon/2}^\text{orb}(\Gamma_R^{(N)}(\mathcal{O}),\|\cdot\|).
    \end{align*}
    Infimizing over $\mathcal{O}$ yields:
    \begin{align*}
        \inf_{E,\delta} \limsup_{N \rightarrow \infty} \frac{1}{N^2} \log K_{F,\epsilon}^\text{orb}( \Gamma_R^{(N,\text{top})}(a;E,\delta),\|\cdot\|) \leq  \inf_{\mathcal{O} \supset L_a}\limsup_{N \rightarrow \infty} \frac{1}{N^2} \log K_{F,\epsilon/2}^\text{orb}(\Gamma_R^{(N)}(\mathcal{O}),\|\cdot\|).
    \end{align*}
    Now, taking the supremum over $\epsilon >0$ and finite subsets $F$ of $I$ completes the proof. 
\end{proof}

\section{The Entropy of Operator Systems}

Recall that an operator system is a subspace $V$ of a $C^*$-algebra which contains the unit, and is closed under the $*$ operation (see \cite{Paulsen} or \cite{BOzawa} for more background on operator systems). In this section, the entropy of operator systems is defined in terms of the $C^*$ algebras they generate. The crucial point is that the entropy, though a priori it depends on the $C^*$-algebra, may be computed on a basis of the operator system.

\medskip

Due to the fact that the natural norm on an operator system is the operator norm, the topological entropy of operator systems is treated first, with the standard $1$-bounded entropy following later in the section. 

\medskip

\noindent \textbf{Proposition 9.1:} Suppose that $\mathscr{A}$ and $\mathscr{B}$ are $C^*$-algebras, and that $\mathscr{A} = C^*(1,S)$ for some subset $S$ of $\mathscr{A}$. Consider the operator system:
\begin{align*}
    V \defeq \text{span} \, (\{1\} \cup \{x,x^*,x^*x,xx^* \, | \, x \in S\}). 
\end{align*}
Suppose that $\Phi: V \rightarrow \mathscr{B}$ is a u.c.p. map for which $\Phi(xx^*)=\Phi(x)\Phi(x)^*$ and $\Phi(x^*x) = \Phi(x)^* \Phi(x)$. Then $\Phi$ is the restriction to $V$ of a $*$-homomorphism $\psi: \mathscr{A} \rightarrow B$. 
\begin{proof}
    By Paulsen 3.6 \cite{Paulsen}, the map $\Phi$ is completely bounded with $1=\|\Phi(1)\| = \|\Phi\| = \|\Phi\|_\text{cb}$. Indeed, $\Phi$ is unital. Since $\|\Phi\|_\text{cb} = 1$, the map $\Phi$ extends to a map $\Tilde{\Phi} : \overline{V} \rightarrow B$ on the closure of $V$ with $\|\Tilde{\Phi}\|_\text{cb} = \|\Tilde{\Phi}\| = 1$. Since $\|\Tilde{\Phi}\|_\text{cb} = \|\Tilde{\Phi}\|$, the extension $\Tilde{\Phi}$ remains u.c.p. Hence, without loss of generality, $V$ may be taken closed.

    \medskip

    Viewing $\mathscr{B}$ as concretely represented in $B(\mathscr{H})$ for some Hilbert space $\mathscr{H}$ via the GNS construction, the Arveson Extension Theorem yields a c.c.p. extension $\psi: \mathscr{A} \rightarrow B(\mathscr{H})$. Since $\{1\} \cup S$ is contained in the multiplicative domain of $\Phi$ by definition and $\{1\} \cup S$ generates a dense $*$-subalgebra, $\psi$ is a $*$-homomorphism. Moreover, on the dense $*$-subalgebra generated by $\{1\} \cup S$, the map $\psi$ is valued in $\mathscr{B}$, so that completeness of $\mathscr{B}$ yields that $\psi$ is valued in $\mathscr{B}$. 
\end{proof}

\medskip

\noindent \textbf{Remark:} A frequently encountered situation in practice is when $S \subset U(\mathscr{A})$, the unitary group of $\mathscr{A}$. In this situation, the hypothesis that $\Phi(xx^*)=\Phi(x)\Phi(x)^*$ and $\Phi(x^*x) = \Phi(x)^* \Phi(x)$ for $x$ in $\{1\} \cup S$ simply states that the images of the members of $\{1\} \cup S$ remain unitaries in $\mathscr{B}$. 

\medskip

Now, to define the entropy of an operator system. Let $V_0$ be an operator system with a spanning set $S$ such that $S^*=S$ and $S$ contains the unit $1$. Intuitively, entropy arises from matrix approximations, and the appropriate maps in this setting are u.c.p. maps. Therefore, setting $B$ to be the norm ultraproduct:
\begin{align*}
    B \defeq \prod_{n \rightarrow \omega} M_{k(n)}(\C),
\end{align*}
the appropriate microstates should arise from u.c.p. maps $\text{UCP}(V_0,B)$. Setting $V$ to be as in the proposition, note the following inclusions:
\begin{align*}
    \text{UCP}(V_0,B) \supset \text{UCP}(V,B) \simeq \text{Hom}(C^*(V),B).
\end{align*}
Indeed, a u.c.p. map from $V$ to $B$ admits a u.c.p. restriction to $V_0$. Moreover, the final identification comes from the above Proposition. Hence, however the entropy of $V_0$ is defined, it must be bounded below by the topological $1$-bounded entropy of $C^*(V)$. 

\medskip

In fact, there is a canonical choice of $C^*(V)$ for which equality makes sense as a definition. Consider $C_\text{max}^*(V_0)$, the $C^*$ algebra generated by $V_0$ and enjoying the universal property that every u.c.p. map $\Phi: V_0 \rightarrow B(\mathscr{H})$ extends to a unique *-homomorphism $C_\text{max}^*(V_0) \rightarrow B(\mathscr{H})$ (see Davidson and Kennedy 4.5 \cite{Ken}). The universal property may be restated as allowing an arbitrary codomain. Indeed, $C^*$-algebras are complete, and the GNS construction allows each to be realized within $B(\mathscr{H})$. Hence,
\begin{align*}
    \text{UCP}(V_0,B) = \text{Hom}(C_\text{max}^*(V_0),B).
\end{align*}
This leads to the following definition.

\medskip

\noindent \textbf{Definition 9.1 (Topological $1$-Bounded Entropy of Operator Systems):} Let $V_0$ be an operator system. Define its topological $1$-bounded entropy $h^\text{top}(V_0)$ to be the topological $1$-bounded entropy of the associated maximal $C^*$ algebra $C_\text{max}^*(V_0)$. 

\medskip

This definition is not a cheap means of extending the above definition for $C^*$-algebras. To demonstrate this, a means of computing $h^\text{top}(V_0)$ on a self-adjoint spanning set $S$ (containing the unit $1$) for $V_0$ is exhibited. In this way, $h^\text{top}(V_0)$ is realized as an \textit{intrinsic} invariant of the operator space $V_0$.

\medskip

To do this, one must introduce the appropriate microstate spaces for $S$ and identify them with topological microstate spaces for $S$ viewed as a tuple in $C_\text{max}^*(V_0)$. 

\medskip

\noindent For this, consider first a u.c.p. embedding $\Phi: V_0 \hookrightarrow B$. Consider $x$ in $S$ and denote $\Phi(x)=(A_x^{(N)})_{N \rightarrow \omega}$. Now, since $\Phi(1)=1$, one has the the norm condition $\|\Phi\|_\text{cb}=1$, and for all natural numbers $m$, constants $\delta>0$, and $(B_y)_{y \in 1 \cup S} \in (M_m(\C))^{1 \cup S}$ there holds:
\begin{align*}
    \bigg\|B_1 \otimes 1 + \sum_{x \in S} B_x \otimes A_x^{(N)}\bigg\| \leq \delta + \bigg\| B_1 \otimes 1 + \sum_{x \in S} B_x \otimes x \bigg\|.
\end{align*}
(for all N in $\omega$). Here, the norm on the left hand side is the operator norm on matrices. On the right hand side, it is the operator norm in matrices over $V_0$. 

\medskip

\noindent \textbf{Definition 9.2:} Fix a finite subset $F$ of the disjoint union:
\begin{align*}
    \bigcup_{m \geq 1} (M_m(\C))^{S}
\end{align*}
and $\delta>0$. Recall that $S$ contains the unit. Define the topological microstates space $\Gamma^{(N,\text{top})}(S;F,\delta)$ to be the space of tuples $(A_x^{(N)})_{x \in S}$ of $N \times N$ matrices over $\C$ such that 
\begin{align*}
    \bigg\|B_1 \otimes 1 + \sum_{x \in S \backslash \{1\}} B_x \otimes A_x^{(N)}\bigg\| \leq \delta + \bigg\| B_1 \otimes 1 + \sum_{x \in S\backslash \{1\}} B_x \otimes x \bigg\|
\end{align*}
for each $(B_y)_{y \in  S}$ belonging to $F$. 

\medskip

Now, it will be demonstrated that the topological entropy of $V_0$ may be computed using these topological microstates spaces. To do this, these microstates spaces must be compared with the topological microstates spaces of $C^*_\text{max}(V_0)$. The following Lemma achieves this.

\medskip

First, some setup is required. Let $V_0$ be an operator system, and let $S$ be a self adjoint spanning set for $V_0$ which contains the unit $1$. Let $\C^*\langle T_y \, | \, y \in S \rangle$ be the noncommutative $*$-algebra on variables labelled by members of $S$. Here, the $*$-operation is given by $T_y^* \defeq T_{y^*}$. 

\medskip

\noindent \textbf{Lemma 9.1:} The following are true.
\begin{enumerate}
    \item For each finite subset $G$ of $\C^*\langle T_y \, | \, y \in S  \rangle$ and $\delta>0$ there exists a finite subset
\begin{align*}
    F' \subset \bigcup_{m \geq 1} M_m(\C)^S
\end{align*}
and $\delta'>0$ such that $\Gamma^{(N,\text{top})}(S,F',\delta') \subset \Gamma^{(N,\text{top})}(S,G,\delta)$ for all sufficiently large $N$. Here, on the right hand side $S$ is viewed as a tuple generating $C^*_\text{max}(V_0)$.
\item On the other hand, for each finite subset
\begin{align*}
    F \subset \bigcup_{m \geq 1} M_m(\C)^S
\end{align*}
and $\delta >0$, there exists a finite subset $G'$ of $\C^*\langle T_y \, | \, y \in S  \rangle$ and $\delta'>0$ such that $\Gamma^{(N,\text{top})}(S,G',\delta') \subset \Gamma^{(N,\text{top})}(S,F,\delta)$ for all sufficiently large $N$. 
\end{enumerate}
\begin{proof}
    First, consider (1). A diagonal argument via contradiction is the approach. Suppose for contradiction that there is a sequence $N_k \rightarrow \infty$ such that for all $\delta'$ and $F'$ as above one has
    \begin{align*}
        \Gamma^{(N_k,\text{top})}(S,F',\delta') \nsubset \Gamma^{(N_k,\text{top})}(S,G,\delta).
    \end{align*}
    for all $k$. Let $\delta_k' \rightarrow 0$, and fix an increasing sequence of finite subsets:
    \begin{align*}
        F_k' \subset \bigcup_{m \geq 1} M_m(\C)^S
    \end{align*}
    such that $[\bigcup_k F_k'] \cap M_m(\C)^S$ is a dense $\Q(i)$-subalgebra for each $m$. For each $k$ let $(A_{k,x})_{x \in S}$ be a member of $\Gamma^{(N_k,\text{top})}(S,F_k',\delta_k')$ which does not belong to $\Gamma^{(N_k,\text{top})}(S,G,\delta)$. Let $\omega$ be a non-principal ultrafilter and define:
    \begin{align*}
        \Phi: V_0 \rightarrow \prod_{k \rightarrow \omega} M_{N_k}(\C)
    \end{align*}
    by $\Phi(x) \defeq (A_{k,x})_{k \rightarrow \omega}$. By the definition of the topological microstates spaces, the fact that $\bigcup_{k} F_k'$ addresses a dense set of matrices, and the choice of $\delta_k'$, one has that $\Phi$ defines a u.c.p. map from $V_0$ into the norm ultraproduct of matrices above. 

    \medskip

    Now, observe that $(A_{k,x})_{x \in S}$ does not belong to $\Gamma^{(N_k,\text{top})}(S,G,\delta)$ for each $k$. Therefore, $\Phi(S)$ is not the image of $S$ under a $*$-homomorphism from $C^*_\text{max}(V_0)$ into the above matrix ultraproduct. This contradicts the universal property of $C^*_\text{max}(V_0)$. 

    \medskip

    A similar argument demonstrates the claim for (2). Only now it is even simpler, because a $*$-homomorphism is evidently a u.c.p. map as a $*$-homomorphism $\psi: \mathscr{A} \rightarrow B$ between $C^*$ algebras $\mathscr{A}$ and $B$ induces further $*$-homomorphisms $\psi \otimes \text{Id}_n : M_n(\mathscr{A}) \rightarrow M_n(B)$. 
\end{proof}

A direct proof of (1) may be given using the linearization trick of Pisier \cite{Lin}.

\begin{proof}
    Fix a finite subset $G$ of $\C^*\langle T_y \, | \, y \in S \rangle$ and $\delta>0$. The claim will be proved under the assumption that $G$ consists of a single $*$-polynomial $P$. Once this is proved, it is evident how to extend the claim to an arbitrary finite subset $G$: simply take the union of the sets $F'$ and the minimum of the $\delta'$ for each $P$ in $G$. Write the members of $S$ as $(x_j)$. Now, by Pisier's linearization \cite{Lin}, there exists a factorization:
    \begin{align*}
        P = P_1 P_2 \dots P_\ell
    \end{align*}
    of $P$ into matrix polynomials of the form:
    \begin{align*}
        P_i(x) = \sum_j a_{j,i} \otimes x_j
    \end{align*}
    (recall that $S$ contains the unit). The $a_{j,i}$ are square matrices. Moreover, for each $\epsilon > 0$, such a factorization exists satisfying:
    \begin{align*}
        \prod_i \|P_i(S)\| \leq (1+\epsilon) \|P(S)\|.
    \end{align*}
    Now, let $\delta'$ and $\epsilon$ both be sufficiently small (what this means will be evident later). Observe that each $(a_{j,i})_{j \geq 1}$ is a member of:
    \begin{align*}
        \bigcup_{m \geq 1} M_m(\C)^S.
    \end{align*}
    Let $F'$ consist of all of the tuples $(a_{j,i})_{j \geq 1}$ as $i$ ranges from $1$ to $\ell$. Now, fix $N$ and suppose that $A=(A_j)_{j \geq 1}$ belongs to $\Gamma^{(N,\text{top})}(S,F',\delta')$. Now, by submultiplicativity of the norm and the definition of the microstates space:
    \begin{align*}
        \|P(A)\| \leq \prod_i \|P_i(A)\| \leq \prod_i (\delta' + \|P_i(S)\|) = O(\delta') + \prod_i \|P_i(S)\| \leq O(\delta') + (1+\epsilon)\|P(S)\| \leq \delta + \|P(S)\|
    \end{align*}
    Here, the last inequality follows from the assumption that $\delta'$ and $\epsilon$ were sufficiently small. Hence, $A$ belongs to $\Gamma^{(N,\text{top})}(S,G,\delta)$. This completes the proof. Once again, (2) is simple even from a constructive viewpoint. The entries of a matrix polynomial are $*$-polynomials.
    
\end{proof}

\medskip

With this Lemma in hand, it will be demonstrated that topological microstates for the operator system $V_0$ may be used to compute the topological entropy of $V_0$. 

\medskip

\noindent \textbf{Theorem 9.1:} Let $V_0$ be an operator system, and $S$ a self-adjoint spanning set for $V_0$ containing the unit. Then
\begin{align}
    h^\text{top}(V_0) = \sup_{\epsilon > 0,E} \inf_{\delta,F} \limsup_{N \rightarrow \infty} \frac{1}{N^2} \log K_{E,\epsilon}^\text{orb}(\Gamma^{(N,\text{top})}(S,F,\delta)).
\end{align}
Here, $E$ ranges over a finite subset of indices indexing $S$, and $F$ and $\delta$ are as above.
\begin{proof}
    Fix a finite subset $G$ of star polynomials $\C^*\langle T_y \, | \,y \in S \rangle$ and $\delta>0$. Now, (1) in the Lemma guarantees that there exists a finite set $F'$ of matrix coefficients and $\delta'>0$ such that for each finite set $E$ of indices:
    \begin{align*}
        \limsup_{N \rightarrow \infty} \frac{1}{N^2}\log K_{E,2\epsilon}^\text{orb}(\Gamma^{(N,\text{top})}(S,F',\delta')) \leq \limsup_{N \rightarrow \infty} \frac{1}{N^2} \log K_{E,\epsilon}^\text{orb}(\Gamma^{(N,\text{top})}(S,G,\delta)).
    \end{align*}
    Note that the $2\epsilon$ is required on the left hand side as the centers of the balls in the packing are required to be within the set to be packed. Now, take the infimum over finite sets $F$ of matrix coefficients and $\delta'>0$. Hence,
    \begin{align*}
        \inf_{F,\delta'}\limsup_{N \rightarrow \infty} \frac{1}{N^2}\log K_{E,2\epsilon}^\text{orb}(\Gamma^{(N,\text{top})}(S,F,\delta')) \leq \limsup_{N \rightarrow \infty} \frac{1}{N^2} \log K_{E,\epsilon}^\text{orb}(\Gamma^{(N,\text{top})}(S,G,\delta)).
    \end{align*}
    Now, take the infimum over $\delta$ and $G$ to obtain:
    \begin{align*}
        \inf_{F,\delta'}\limsup_{N \rightarrow \infty} \frac{1}{N^2}\log K_{E,2\epsilon}^\text{orb}(\Gamma^{(N,\text{top})}(S,F,\delta')) \leq \inf_{G,\delta}\limsup_{N \rightarrow \infty} \frac{1}{N^2} \log K_{E,\epsilon}^\text{orb}(\Gamma^{(N,\text{top})}(S,G,\delta)).
    \end{align*}
    Finally, taking the supremum over $\epsilon$ yields that the quantity on the RHS in (1) does not exceed $h(V_0)$. The reverse inequality is similar, using (2) in the Lemma instead.
    
\end{proof}

\medskip

Now that the topological entropy was computed on a spanning set of an operator system, it is time to demonstrate the analogous claim for $1$-bounded entropy. To do this, it is necessary to define the $1$-bounded entropy of an operator system. In this setting, embeddings into the tracial ultraproduct are under consideration:
\begin{align*}
    \mathscr{M} \defeq \prod_{n \rightarrow \omega} (M_{k(n)}(\C),\text{tr}_{k(n)}). 
\end{align*}
Let $V$ be an operator system, and $S$ a self-adjoint spanning set for $S$ which contains the unit. As in the case of topological $1$-bounded entropy, the universal property yields $\text{Hom}(C^*_\text{max}(V),\mathscr{M}) \simeq \text{UCP}(V,\mathscr{M})$. Recall also that a c.b. map $\Phi:V \rightarrow \mathscr{M}$ is u.c.p. provided $\|\Phi\|_\text{cb} \leq 1$, the map $\Phi$ commutes with adjoints, and $\Phi$ is unital. 

\medskip

\noindent \textbf{Definition 9.3 ($1$-Bounded Entropy of Operator Systems):} Let $V$ be an operator system. Define its $1$-bounded entropy $h(V)$ to be the $1$-bounded entropy of the tracial completion of $C^*_\text{max}(V)$ with respect to all of its traces. 

\medskip

Now, it must be demonstrated again that this notion is intrinsic to $V$. As above, the definition of microstates centers around the classification of tuples $(a_s)_{s \in S}$ belonging to $\mathscr{M}^S$ for which there exists a u.c.p. map $\Phi: V \rightarrow \mathscr{M}$ such that $\Phi(s) = a_s$ for each $s$ in $S$. Once again, this translates to the norm condition:
\begin{align*}
    \bigg\|B_1 \otimes 1 + \sum_{y \in S \backslash \{1\}} B_y \otimes a_y \bigg\| \leq \bigg\| B_1\otimes 1 + \sum_{y \in S \backslash \{1\}} B_y \otimes y \bigg\|.
\end{align*}
for all tuples $(B_y)$ belonging to
\begin{align*}
    \bigcup_{m \geq 1} M_m(\C)^S.
\end{align*}
On the left hand side, the norm is that on matrices over $\mathscr{M}$, whereas the norm on the right hand side is evaluated on matrices over $V$. In order to define microstates, one must fix lifts $(A_{N,s})_{N \rightarrow \omega}$ for each of the $a_s$. A notable difference between the $1$-bounded entropy and the topological $1$-bounded entropy is the necessity of a norm cutoff. Henceforth, suppose:
\begin{align*}
    \sup_{N} \|A_{N,s}\| \leq R
\end{align*}
for a fixed $R$. In order to relate the operator norm condition to the trace, simply observe that for each member $a$ of $\mathscr{M}$ there holds:
\begin{align*}
    \|a\| = \sup_{1 \leq q < \infty} \|a\|_{L^q(\text{tr}_\omega)}.
\end{align*}
Furthermore, since continuous functional calculus commutes with taking the ultraproduct (i.e. $f(|a|) = (f(|A_N|))_{N \rightarrow \omega}$ for $(A_N)_{N \rightarrow \omega}$ a lift of $a$), one has for each $1\leq q<\infty$
\begin{align*}
    \|a\|_{L^q(\text{tr}_\omega)} = \lim_{N \rightarrow \omega} \|A_N\|_{L^q(\text{tr}_N)}.
\end{align*}
Here, $f(\cdot) = \text{tr}(|\cdot|^q)$. Putting these pieces together, the above inequality in the norms translates into the following condition: for all finite subsets:
\begin{align*}
    F \subset \bigcup_{m \geq 1} M_m(\C)^S,
\end{align*}
constants $\delta>0$, and $1\leq q<\infty$, there exists an $\omega$ large set $\Omega$ such that for all $N$ belonging to $\Omega$,
\begin{align*}
    \bigg\| B_1 \otimes 1 + \sum_{y \in S \backslash \{1\}} B_y \otimes A_{N,y} \bigg\|_{L^q(\text{tr}_N)} \leq \delta + \bigg\| B_1\otimes 1 + \sum_{y \in S \backslash \{1\}} B_y \otimes y \bigg\|.
\end{align*}
Hence, the following definition.

\medskip

\noindent \textbf{Definition 9.4:} Let $R>0$ be an operator norm cutoff. Fix a finite subset 
\begin{align*}
    F \subset \bigcup_{m \geq 1} M_m(\C)^S,
\end{align*}
a constant $\delta>0$, and $1 \leq q <\infty$. Now, for each natural number $N$, define the microstates space $\Gamma_R^{(N)}(S,F,\delta,q)$ to consist of those tuples $(A_{N,y})_{y \in S}$ of $N \times N$ matrices of operator norm not exceeding $R$ such that 
\begin{align*}
    \bigg\| B_1 \otimes 1 + \sum_{y \in S \backslash \{1\}} B_y \otimes A_{N,y} \bigg\|_{L^q(\text{tr}_N)} \leq \delta + \bigg\| B_1\otimes 1 + \sum_{y \in S \backslash \{1\}} B_y \otimes y \bigg\|
\end{align*}

for all tuples $(B_y)_{y \in S}$ belonging to $F$.

\medskip

As in the case of topological $1$-bounded entropy, it will be demonstrated that these microstates spaces may be used to compute $h(V)$. The analogous Lemma is as follows.

\medskip

\noindent \textbf{Lemma 9.2:} Fix a norm cutoff $R>0$ such that each of the members of $S$ is bounded by $R$ in norm. Recall that the tuple $S$ (viewed within $C^*_\text{max}(V)$) defines a compact, convex subset $L_S$ of the space $\Sigma_{R,S}$ of $R$-bounded laws. The following statements hold.
\begin{enumerate}
    \item Fix a weak-$*$ neighborhood $\mathcal{O}$ of $L_S$. There exists $\delta>0$, a finite subset:
    \begin{align*}
        F \subset \bigcup_{m \geq 1} M_m(\C)^S,
    \end{align*}
    and $1 \leq q < \infty$ so that for all $N$ sufficiently large:
    \begin{align*}
        \Gamma_R^{(N)}(S,F,\delta,q) \subset \Gamma_R^{(N)}(\mathcal{O}).
    \end{align*}

    \item On the other hand, for each $\delta > 0$, each finite subset $F$ as in (1), and each $1 \leq q < \infty$, there exists a weak-$*$ neighborhood $\mathcal{O}$ of $L_S$ such that for all $N$ sufficiently large:
    \begin{align*}
        \Gamma_R^{(N)}(\mathcal{O}) \subset \Gamma_R^{(N)}(S,F,\delta,q).
    \end{align*}
\end{enumerate}
\begin{proof}
    First, a proof of (1) is given. It proceeds via a diagonal argument and contradiction. Suppose for contradiction that there exists $N_k \rightarrow \infty$ such that for all $\delta,F,q$ as above,
    \begin{align*}
        \Gamma_R^{(N_k)}(S,F,\delta,q) \nsubset \Gamma_R^{(N_k)}(\mathcal{O})
    \end{align*}
    for all $k$. Let $\delta_k \rightarrow 0$, $q_k \rightarrow \infty$, and let $F_k$ be an increasing sequence of finite subsets:
    \begin{align*}
        F_k \subset \bigcup_{m \geq 1} M_m(\C)^S
    \end{align*}
    for which $[\bigcup_k F_k] \cap M_m(\C)^S$ is a dense $\Q(i)$-subalgebra for each $m$. For each $k$, let $(A_{k,x})_{x \in S}$ be a member of $\Gamma_R^{(N_k)}(S,F_k,\delta_k,q_k)$ which does not belong to $\Gamma_R^{(N_k)}(\mathcal{O})$. Define a map:
    \begin{align*}
        \Phi: V \rightarrow \prod_{k \rightarrow \omega} (M_{N_k}(\C),\text{Tr}_{N_k}) 
    \end{align*}
    by $\Phi(x) \defeq (A_{k,x})_{k \rightarrow \omega}$. The map $\Phi$ is well-defined and u.c.p. by the choices of $F_k$, $\delta_k$, and $q_k$ made above. Indeed, this follows by the definition of the microstates spaces. By the universal property of $C^*_\text{max}(V)$, the map $\Phi$ extends to a $*$-homomorphism \cite{Ken}:
    \begin{align*}
        C^*_\text{max}(V) \rightarrow \prod_{k \rightarrow \omega} (M_{N_k}(\C),\text{Tr}_{N_k}).
    \end{align*}
    Pulling back the trace $\text{Tr}_\omega \defeq \lim_{k \rightarrow \omega} \text{Tr}_{N_k}$ on the ultraproduct yields a trace $\tau$ on $C^*_\text{max}(V)$ such that the law of $S$ with respect to $\tau$ coincides with the law of $[(A_{k,x})_{k \rightarrow \infty}]_{x \in S}$ with respect to $\text{Tr}_\omega$. Therefore, the laws of $(A_{k,x})_{x \in S}$ (denote them $\ell_k$) with respect to $\text{Tr}_{N_k}$ converge weakly to the law of $S$ with respect to $\tau$. Thus, there is $k_0$ such that $\ell_k$ belongs to $\mathcal{O}$ for all $k>k_0$. However, this implies that for large enough $k$, one has $(A_{k,x})_{x \in S}$ belongs to $\Gamma_R^{(N_k)}(\mathcal{O})$, a contradiction. 

    \medskip

    The converse proceeds by a similar contradiction. Choosing microstates arising from a shrinking sequence of weak-$*$ neighborhoods, one obtains a trace preserving embedding $C^*_\text{max}(V)$ into a matrix ultraproduct. The map is u.c.p. and lifts of the images of members of $S$ form microstates for $V$ in the above defined sense. 
\end{proof}

\medskip

\noindent \textbf{Theorem 9.2:} Fix a norm cutoff $R>0$. Let $V$ be an operator system and $S$ a self-adjoint spanning set for $V$ containing the unit such that $\|x\| \leq R$ for each member $x$ of $S$. Then
\begin{align}
    h(V) = \sup_{\epsilon > 0,E} \inf_{F,\delta,q} \limsup_{N \rightarrow \infty} \frac{1}{N^2} \log K_{E,\epsilon}^\text{orb}(\Gamma_R^{(N)}(S,F,\delta,q))
\end{align}
where $E$ ranges over finite subsets of the set indexing $S$. 
\begin{proof}
    Fix a weak-$*$ neighborhood $\mathcal{O}$ of $L_S$. By the Lemma above, there exists $F$, $\delta$, and $q$ as above such that
    \begin{align*}
        \limsup_{N \rightarrow \infty} \frac{1}{N^2}\log K_{E,2\epsilon}^\text{orb}(\Gamma_R^{(N)}(S,F,\delta,q)) \leq \limsup_{N \rightarrow \infty}\frac{1}{N^2}\log K_{E,\epsilon}^\text{orb}(\Gamma_R^{(N)}(\mathcal{O})).
    \end{align*}
   Again, recall that due to the restriction on the centers in the packing, $2\epsilon$ is necessary. Now, take the infimum over $F$, $\delta$, and $q$ to obtain:
    \begin{align*}
       \inf_{F,\delta,q} \limsup_{N \rightarrow \infty} \frac{1}{N^2}\log K_{E,2\epsilon}^\text{orb}(\Gamma_R^{(N)}(S,F,\delta,q)) \leq \limsup_{N \rightarrow \infty}\frac{1}{N^2}\log K_{E,\epsilon}^\text{orb}(\Gamma_R^{(N)}(\mathcal{O})).
    \end{align*}
    Now, infimize over the weak-$*$ neighborhoods $\mathcal{O}$:
    \begin{align*}
        \inf_{F,\delta,q} \limsup_{N \rightarrow \infty} \frac{1}{N^2}\log K_{E,2\epsilon}^\text{orb}(\Gamma_R^{(N)}(S,F,\delta,q)) \leq \inf_{\mathcal{O} \supset L_S}\limsup_{N \rightarrow \infty}\frac{1}{N^2}\log K_{E,\epsilon}^\text{orb}(\Gamma_R^{(N)}(\mathcal{O})).
    \end{align*}
    Finally, take the supremum over $\epsilon>0$ and $E$ to obtain that the RHS in (2) does not exceed $h(V)$. For the reverse inequality, the proof is similar, utilizing part 2 in the above Lemma.
    
\end{proof}

Now that the tools have been developed, some applications are laid out in the remainder of the paper. The focus is on computations with $C^*$-algebras, although in principal these results also hold implications for operator systems. While most results involve the $1$-bounded entropy, keep in mind the inequality $h^\text{top} \leq h$ for applications to topological entropy. Moreover, it is worth mentioning that the focus here is on studying $C^*$-algebras via their tracial completions. There are specializations of von Neumann algebraic properties (such as the McDuff property) in the setting of tracially complete $C^*$ algebras which are not explored here \cite{TC}.

\section{Some Vanishing Conditions}

In this section, some computations of the $1$-bounded entropy are carried out, and some implications thereof are developed. 

\medskip

Recall that the $1$-bounded entropy of a hyperfinite von Neumann algebra vanishes. Hence, a similar vanishing result for the entropy of $C^*$ algebras is suggested in the nuclear case. The proof illustrates the general theme that the variational principle allows for reductions to known results regarding von Neumann algebras. 

\medskip

\noindent \textbf{Theorem 10.1:} Let $\mathscr{A}$ be a nuclear $C^*$ algebra. Then $h(\mathscr{A})=0$.
\begin{proof}
    Applying the variational principle, 
    \begin{align*}
        h(\mathscr{A}) = \sup_{\tau \in \mathscr{T}(\mathscr{A})} h(M_\tau),
    \end{align*}
    where $M_\tau$ denotes the von Neumann algebra obtained via the GNS construction associated to the trace $\tau$. Since $\mathscr{A}$ is nuclear, $M_\tau$ is hyperfinite for each trace $\tau$ on $\mathscr{A}$. Therefore, $h(M_\tau)=0$ for each such trace $\tau$ (this is an exercise from the orbital packing definition in Hayes section 2 \cite{Reg}), so that $h(\mathscr{A})=0$.
\end{proof}

Observe that via the above comparison result, one also obtains $h^\text{top}(\mathscr{A}) \leq 0$. As a brief example, recall that the rotation algebras $A_\theta$ are nuclear, so that $h(A_\theta)=0$. 

\medskip


A few more related conditions for vanishing of $h$ are considered here. First, recall that a $C^*$-algebra $\mathscr{A}$ is \textit{Jiang-Su stable} provided $\mathscr{A} \otimes_\text{min} \mathscr{Z} \cong \mathscr{A}$, where $\mathscr{Z}$ is the Jiang-Su algebra: a cornerstone in the Elliot classification program (see for instance \cite{Classification}). 

\medskip

For the present purpose, the reader need only know that $\mathscr{Z}$ supports a unique trace, and that the GNS representation of $\mathscr{Z}$ with respect to this trace generates a copy of the hyperfinite II$_1$ factor $R$. Indeed, this trace is faithful due to the fact that $\mathscr{Z}$ is simple. In this way, the property of Jiang-Su stability is the analogue of the McDuff property for II$_1$ factors in the setting of $C^*$-algebras. The entropy of a McDuff II$_1$ factor vanishes by the primeness results in \cite{Reg} (a McDuff factor tensorially absorbs $R$), and in this setting the analogous Corollary is obtained. 

\medskip

\noindent \textbf{Corollary 10.1:} Suppose $\mathscr{A}$ is a Jiang-Su stable $C^*$-algebra. Then $h(\mathscr{A}) \leq 0$.
\begin{proof}
    Identify $\mathscr{A}$ with the minimal tensor product $\mathscr{A} \otimes_\text{min} \mathscr{Z}$ and fix a trace $\tau$ on $\mathscr{A} \otimes_\text{min} \mathscr{Z}$. The restriction of $\tau$ to the second tensorand must be the unique trace on $\mathscr{Z}$. Hence, the GNS von Neumann algebra is:
    \begin{align*}
    M_\tau \defeq \pi_{\tau|_\mathscr{A}}(\mathscr{A})'' \otimes R.
    \end{align*}
    If $\pi_{\tau|_\mathscr{A}}(\mathscr{A})''$ has diffuse center, then so does $M_\tau$ so that $h(M_\tau) \leq 0$. Therefore, suppose without loss of generality that $\pi_{\tau|_\mathscr{A}}(\mathscr{A})''$ has atomic center. Indeed, $M_\tau$ splits as a direct sum of von Neumann algebras, one with diffuse center and the other with atomic center; the $1$-bounded entropy of a tracial von Neumann algebra satisfies a weighted subadditivity with respect to direct sums \cite{Reg} A.13. Since tensor products distribute over direct sums:
    \begin{align*}
        M_\tau \cong \bigoplus_{i \in I} (z_i\pi_{\tau|_\mathscr{A}}(\mathscr{A})'' \otimes R),
    \end{align*}
    where the $z_i$ are pairwise orthogonal central projections in $\pi_{\tau|_\mathscr{A}}(\mathscr{A})''$ and their corresponding corners $z_i\pi_{\tau|_\mathscr{A}}(\mathscr{A})''$ are factors. Each of the summands is McDuff, as it is a tensor product with $R$, the tensor product operation is associative, and $R \otimes R \cong R$. Indeed, it is evident that tensorial absorption of $R$ implies the McDuff property as $R$ itself has highly nontrivial central sequences. Recall that the $1$-bounded entropy of a McDuff factor is at most $0$ (this follows from primeness results in \cite{Reg}). Therefore, the weighted subadditivity of the $1$-bounded entropy under direct sums \cite{Reg} A.13 yields the following:
    \begin{align*}
        h(M_\tau) \leq \sum_{i \in I} \tau(z_i)^2h(z_i\pi_{\tau|_\mathscr{A}}(\mathscr{A})'' \otimes R) =0.
    \end{align*}
    Appealing to the variational principle, $h(\mathscr{A})\leq 0$. 
    
\end{proof}

The above argument illustrates the relevance of the unique trace property in computations of the $1$-bounded entropy for $C^*$ algebras. However, one also sees immediately that the argument may be generalized to $C^*$ algebras $\mathscr{A}$ which enjoy the following weakening of uniform property $\Gamma$: each II$_1$ factor representation of $\mathscr{A}$ has property $\Gamma$. 

\medskip

\noindent \textbf{Corollary 10.2:} Suppose that $\mathscr{A}$ is a $C^*$ algebra such that each II$_1$ factor representation of $\mathscr{A}$ has property $\Gamma$. Then $h(\mathscr{A}) \leq 0$.
\begin{proof}
    Fix a trace $\tau$ on $\mathscr{A}$ and consider its associated GNS von Neumann algebra $M_\tau$. If $M_\tau$ has diffuse center, then automatically $h(M_\tau) \leq 0$. Suppose without loss of generality that $M_\tau$ has atomic center. Indeed, $M_\tau$ is a direct sum of von Neumann algebras, one of which has diffuse center and one of which has atomic center; the $1$-bounded entropy of a tracial von Neumann algebra satisfies a weighted subadditivity with respect to direct sums \cite{Reg} A.13. Then there is a family $(z_i)_{i \in I}$ of pairwise orthogonal central projections such that:
    \begin{align*}
        M_\tau \cong \bigoplus_{i \in I} z_iM_\tau,
    \end{align*}
    so that $z_iM_\tau$ is a factor for each $i$. Since the trace restricts to each of these corners (albeit scaled such that it remains a state), one has that $z_i M_\tau$ is a finite factor for each $i$. Hence, if $z_i M_\tau$ is diffuse, it is a II$_1$ factor.

    \medskip
    
    In the case in which $z_i M_\tau$ is diffuse, the hypothesis yields that $z_i M_\tau$ has property $\Gamma$. This yields $h(z_iM_{\tau})\leq 0$ (see for instance \cite{Absorption} Example 4). On the other hand, if $z_i M_\tau$ is a finite dimensional factor, it is isomorphic to a matrix algebra. Any two embeddings of a matrix algebra into a II$_1$ factor are conjugate (simply choose matrix units and use comparison of projections to produce partial isometries relating them), such that $h(z_iM_\tau)=0$ in this case. Hence, by the weighted subadditivity property: 
    \begin{align*}
        h(M_\tau) \leq \sum_{i \in I} \tau(z_i)^2 h(z_i M_\tau) \leq 0.
    \end{align*}
    Now, apply the variational principle to obtain $h(\mathscr{A}) \leq 0$. 
\end{proof}

\medskip

When discussing the vanishing of the entropy, the proof of the Peterson-Thom conjecture yields some natural Corollaries. The conjecture was proved to be equivalent to the following $1$-Bounded Entropy Conjecture by Hayes \cite{HayesThom}: if $Q$ is a diffuse subalgebra of a free group factor $L(\F_r)$ such that $h(Q)=0$, then $Q$ must be amenable \cite{PT}.

\medskip

This result has a recent and exciting history. In groundbreaking work of Hayes \cite{HayesThom}, a random matrix criterion implying the conjecture was posited. Subsequently, this criterion was verified using strong convergence methods (see \cite{BelCap}, \cite{BorCol}, \cite{MagSal}, \cite{Parraud}, and \cite{NewApproach}). The following Corollary provides a $C^*$-algebraic version.

\medskip

\noindent \textbf{Corollary:} Suppose that $\mathscr{A}$ is a $C^*$-subalgebra of $C_\lambda^*(\F_r)$ for some $r \geq 2$ such that $h(\mathscr{A})=0$. Moreover, suppose that $\pi_\tau(\mathscr{A})''$ is diffuse, where $\tau$ is the standard trace on $C_\lambda^*(\F_r)$. Then the restriction of $\tau$ to $\mathscr{A}$ is an amenable trace.
\begin{proof}
    Applying the variational principle:
    \begin{align*}
        h(\pi_\tau(\mathscr{A})'') = 0.
    \end{align*}
    Now, observe that $\pi_\tau(\mathscr{A})'' \subset L(\F_r)$. Applying the $1$-Bounded Entropy Conjecture yields that $\pi_\tau(\mathscr{A})''$ is amenable. Moreover, since the trace $\tau$ is faithful, $\pi_\tau(\mathscr{A})$ is an isometric copy of $\mathscr{A}$ which generates $\pi_\tau(\mathscr{A})''$. Since $Q$ is amenable, there is a conditional expectation $B(L^2(Q)) \rightarrow Q$. Hence, applying \cite{BOzawa} 6.2.7 (see also the footnote contained there) completes the proof. 
\end{proof}

\medskip

As a concrete example, recall that a countable, discrete group $\Gamma$ is amenable if and only if $C_\lambda^*(\Gamma)$ admits an amenable trace \cite{BOzawa}.

\medskip

The condition that $\pi_\tau(\mathscr{A})''$ is diffuse cannot be disposed of. It is ultimately essential in most applications of the variational principle framework. This is the point which will be elaborated in the following section.

\section{The Proper Notion of a Diffuse Subalgebra}

Regarding behavior of the $1$-bounded entropy for von Neumann algebras, an indispensable hypothesis is the property of being diffuse \cite{Reg}. That is, the absence of minimal projections. 

\medskip

Given that a $C^*$ algebra may possess no nontrivial projections, this hypothesis does not carry over naively to the present setting, and some care needs to be taken regarding the correct generalization in the context of $C^*$ algebras. This warrants the discussion of this section.

\medskip

The primary difficulty in generalizing this notion is that one needs control over all of the GNS representations $\pi_\tau$ coming from traces on a $C^*$-algebra. Consider the following extreme example, which illustrates what might go wrong when there is no control over GNS representations.

\medskip

\noindent \textbf{Example 11.1:} Consider the full group $C^*$-algebra $C^*(\F_2 \times \F_2) \cong C^*(\F_2) \otimes_\text{max} C^*(\F_2)$. Recall that $\F_2 \times \F_2$ is a right angled Artin group on the following connected graph:

\begin{center}
\begin{figure}[h!]
    \centering
    \tikzset{every picture/.style={line width=0.75pt}} 

\begin{tikzpicture}[x=0.75pt,y=0.75pt,yscale=-1,xscale=1]

\draw    (297.25,83) -- (345.25,131) ;
\draw    (297.25,83) -- (346.25,82) ;
\draw    (346.25,82) -- (299.25,131) ;
\draw    (299.25,131) -- (345.25,131) ;
\draw  [fill={rgb, 255:red, 0; green, 0; blue, 0 }  ,fill opacity=1 ] (294.13,83) .. controls (294.13,81.34) and (295.52,80) .. (297.25,80) .. controls (298.98,80) and (300.38,81.34) .. (300.38,83) .. controls (300.38,84.66) and (298.98,86) .. (297.25,86) .. controls (295.52,86) and (294.13,84.66) .. (294.13,83) -- cycle ;
\draw  [fill={rgb, 255:red, 0; green, 0; blue, 0 }  ,fill opacity=1 ] (299.25,131) .. controls (299.25,129.34) and (300.65,128) .. (302.38,128) .. controls (304.1,128) and (305.5,129.34) .. (305.5,131) .. controls (305.5,132.66) and (304.1,134) .. (302.38,134) .. controls (300.65,134) and (299.25,132.66) .. (299.25,131) -- cycle ;
\draw  [fill={rgb, 255:red, 0; green, 0; blue, 0 }  ,fill opacity=1 ] (342.13,131) .. controls (342.13,129.34) and (343.52,128) .. (345.25,128) .. controls (346.98,128) and (348.38,129.34) .. (348.38,131) .. controls (348.38,132.66) and (346.98,134) .. (345.25,134) .. controls (343.52,134) and (342.13,132.66) .. (342.13,131) -- cycle ;
\draw  [fill={rgb, 255:red, 0; green, 0; blue, 0 }  ,fill opacity=1 ] (343.13,82) .. controls (343.13,80.34) and (344.52,79) .. (346.25,79) .. controls (347.98,79) and (349.38,80.34) .. (349.38,82) .. controls (349.38,83.66) and (347.98,85) .. (346.25,85) .. controls (344.52,85) and (343.13,83.66) .. (343.13,82) -- cycle ;

\end{tikzpicture}
    \caption{The graph corresponding to the RAAG $\F_2 \times \F_2$.}
    \label{fig:placeholder}
\end{figure}
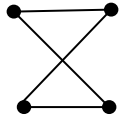
\end{center}
As such, it is generated by sequentially commuting elements. Despite the fact that this $C^*$-algebra is a tensor product (albeit a maximal one), and has a "quasi-regular" (sequential commutation) copy of the circle algebra $C(S^1)$ contained within it, its entropy is still infinite for the simple reason that it surjects:
\begin{align*}
    C^*(\F_2 \times \F_2) \rightarrow C^*_\lambda(\F_2)
\end{align*}
onto the reduced group $C^*$-algebra of $\F_2$. Indeed, the $1$-bounded entropy of $L(\F_2)$ is infinite (\cite{Voic2}, \cite{Jung}, and \cite{Reg}), so applying the variational principle completes the computation.

\medskip

This example illustrates that when investigating inclusions $\mathscr{B} \subset \mathscr{A}$ of $C^*$ algebras involving products and regularity, one must require that for each trace $\tau$ on $\mathscr{A}$, the subalgebra $\mathscr{B}$ generates a diffuse von Neumann algebra $\pi_\tau(\mathscr{B})''$. This fails in the present example for each choice of relevant subalgebra. For instance, the relevant copy of the circle algebra $C(S^1)$ fell within the kernel of the above surjection. Or, if one of the tensorands $C^*(\F_2)$ is considered, each of these hosts many finite dimensional representations compatible with traces on the entire algebra.

\medskip

In view of the variational principle, the above mentioned generalization of the diffuse hypothesis is also quite practical, as will be seen below. Before proceeding to more specific/delicate settings, it is useful to mention a few generalities which guarantee diffuse GNS von Neumann algebras.

\medskip

\noindent \textbf{Proposition 11.1:} Fix a $C^*$-algebra $\mathscr{A}$. Then the GNS von Neumann algebra $\pi_\tau(\mathscr{A})''$ is diffuse for each trace $\tau$ on $\mathscr{A}$ if and only if for all traces $\tau$ on $\mathscr{A}$ one has $\pi_\tau(\mathscr{A})$ infinite-dimensional.
\begin{proof}
    First, suppose that there exists a trace $\tau$ on $\mathscr{A}$ such that $\pi_\tau(\mathscr{A})$ is finite dimensional. A finite-dimensional $C^*$ algebra is already a von Neumann algebra, such that $\pi_\tau(\mathscr{A})''$ is atomic.

    \medskip

    Conversely, suppose that there exists a trace $\tau$ on $\mathscr{A}$ such that $M \defeq \pi_\tau(\mathscr{A})''$ is not diffuse. Hence, there exists a minimal projection $p$ belonging to $M$. Observe that since $p$ is a minimal projection, it is abelian: $pMp \cong \C p$. Pass to the central support $z(p)$ and consider the corner $z(p)M$. 

    \medskip

    It will quickly be demonstrated that $z(p)M$ is type I. If $z$ is a nonzero central projection, then $z$ and $z(p)$ are clearly not orthogonal. Hence, there exists a subprojection of $p$ and a subprojection $q$ of $z$ which are equivalent. But $p$ is minimal, so that $p \sim q$. Hence, $q$ is abelian. Since $z$ was arbitrary, $z(p)M$ has type I. 

    \medskip

    Now, note that $z(p)M$ carries the trace $\tau_p$ defined by $\tau_p(x) \defeq \tau(z(p)x)/\tau(z(p))$. Hence, $z(p)M$ is of finite type I. Moreover, it is readily demonstrated from the minimality of $p$ that $z(p)$ is minimal in $Z(M)$, so that $z(p)M$ is a factor. It follows that $z(p)M$ is a matrix algebra. Now, the formula $x \mapsto \tau_p(z(p)\pi_\tau(x))$ defines a trace on $\mathscr{A}$ such that the GNS completion is $z(p)M$, which is finite-dimensional.
\end{proof}

\medskip

\noindent \textbf{Proposition 11.2:} Suppose $G$ is a countable, discrete, non-amenable group. Then the reduced group $C^*$-algebra $C_\lambda^*(G)$ is such that the GNS von Neumann algebra $\pi_\tau(C_\lambda^*(G))''$ is diffuse for each trace $\tau$ on $C_\lambda^*(G)$. 
\begin{proof}
    Examining the proof of the above Proposition, the existence of a trace $\tau$ on $C_\lambda^*(G)$ such that $\pi_\tau(C_\lambda^*(G))''$ is not diffuse implies the existence of a finite-dimensional representation $\pi_{\tau}: C_\lambda^*(G) \rightarrow M_n(\C)$. This induces a representation $\pi$ of the full group $C^*$ algebra of $G$ via the composition with the canonical map from $C^*(G)$ onto $C_\lambda^*(G)$:
    \begin{align*}
        C^*(G) \rightarrow C_\lambda^*(G) \rightarrow M_n(\C).
    \end{align*}
    Now, recall that representations of $C^*(G)$ are in one to one correspondence with unitary representations of $G$. Hence, there is a representation $\rho :G \rightarrow GL_n(\C)$ which extends to the above representation of $C^*(G)$. 

    \medskip

    Observe that since $\pi$ factors through the reduced group $C^*$-algebra, $\rho \prec \lambda$ where $\lambda$ is the left regular representation. By Fell's absorption principle, it follows that $\rho \otimes \rho^* \prec \lambda$ as well. However, since $\rho$ is a finite-dimensional representation, $\rho \otimes \rho^*$ contains a trivial subrepresentation. Since the trivial representation of $G$ is weakly contained in the left regular representation, $G$ is amenable \cite{BOzawa}.
\end{proof}

\medskip

A more abstract condition which guarantees the diffuse hypothesis for all GNS von Neumann algebras corresponding to all traces is as follows.

\medskip

\noindent \textbf{Definition 11.1 (Nowhere-Scattered $C^*$-Algebras \cite{Thiel}):} Fix a $C^*$-algebra $\mathscr{A}$. The algebra $\mathscr{A}$ is said to be nowhere-scattered provided for each extreme point $\varphi$ of the space of states $S(\mathscr{A})$ (I.e. a pure state) there holds:
\begin{align*}
    \pi_\varphi(\mathscr{A}) \cap K(L^2(\varphi)) = \{0\}.
\end{align*}
That is, the GNS $C^*$-algebra $\pi_\varphi$ intersects the compacts trivially. 

\medskip

This property guarantees that each GNS von Neumann algebra $\pi_\tau(\mathscr{A})''$ corresponding to a trace $\tau$ on $\mathscr{A}$ is diffuse. Indeed, this condition is an analogue of the type II condition for factors. A finite von Neumann algebra is nowhere scattered as a $C^*$-algebra if and only if it is of type II \cite{Thiel}.

\section{Regular Subalgebras}

Now that the proper notion of a diffuse subalgebra has been introduced, a discussion of regularity may be pursued. A couple of results of a similar nature to those in work of Hayes \cite{Reg} in the case of tracial von Neumann algebras are obtained. First, crossed products are treated.

\medskip

\noindent \textbf{Theorem 12.1:} Suppose that $\mathscr{A}$ is isomorphic to the reduced crossed product $\mathscr{B} \rtimes_{\alpha,r} \Gamma$ corresponding to the action of a discrete group $\Gamma$ on a $C^*$ algebra $\mathscr{B}$. Suppose additionally, that $\pi_\tau(\mathscr{B})''$ is diffuse for all traces $\tau$ on $\mathscr{A}$. Then:
\begin{align*}
    h(\mathscr{A}) \leq h(\mathscr{B}).
\end{align*}
\begin{proof}
    For each trace $\tau$, the GNS von Neumann algebra:
    \begin{align*}
        \pi_\tau(\mathscr{A})'' \cong \pi_\tau(\mathscr{B})'' \rtimes \Gamma.
    \end{align*}
    Hence, for each trace $\tau$ one has that $\pi_\tau(\mathscr{B})''$ is a diffuse, regular subalgebra of $\pi_\tau(\mathscr{A})''$. Hence, applying Hayes \cite{Reg} and the variational principle:
    \begin{align*}
        h(\mathscr{A}) = \sup_{\tau \in \mathscr{T}(\mathscr{A})} h(\pi_\tau(\mathscr{A})'') \leq \sup_{\tau \in \mathscr{T}(\mathscr{A})} h(\pi_\tau(\mathscr{B})'') \leq \sup_{\tau \in \mathscr{T}(\mathscr{B})} h(\pi_\tau(\mathscr{B})'') = h(\mathscr{B}).
    \end{align*}
    Indeed, each trace on $\mathscr{A}$ restricts to a trace on $\mathscr{B}$. 
\end{proof}

Of particular interest is the application of this result to classical $C^*$ dynamical systems. That is, the case in which the subalgebra $\mathscr{B}$ is $C(X)$ for $X$ a compact metrizable space. In this situation, the condition on diffuse GNS representations is particularly concrete.

\medskip

\noindent \textbf{Corollary 12.1:} Let $X$ be a compact, metrizable topological space and let $\Gamma$ be a countable discrete group with an action $\alpha: \Gamma \curvearrowright X$ on $X$ by homeomorphisms. Observe that $\alpha$ induces an action of $\Gamma$ on $C(X)$. Denote it $\beta$. Suppose additionally that each ergodic probability measure $\mu$ on $X$ is atomless. Then:
\begin{align*}
    h(C(X) \rtimes_{\beta,r} \Gamma) \leq 0. 
\end{align*}
\begin{proof}
    In order to apply the above Theorem, it suffices to demonstrate that the above condition on Ergodic measures guarantees that the GNS von Neumann algebras $\pi_\tau(C(X))''$ are diffuse for all traces $\tau$ on the reduced crossed product. Indeed, the entropy of $C(X)$ is $0$ as it is a nuclear algebra.

    \medskip
    
    Recall that traces $\tau$ on the reduced crossed product correspond directly to $\Gamma$ invariant probability measures $\nu_\tau$ on $X$ in such a way that:
    \begin{align*}
       & \pi_\tau(C(X) \rtimes_{\beta,r} \Gamma)'' \cong L^\infty(X,\nu_\tau) \rtimes \Gamma \\
       & \pi_\tau(C(X))'' \cong L^\infty(X,\nu_\tau).
    \end{align*}
    Observe that if $\nu_\tau$ is atomless, then $L^\infty(X,\nu_\tau)$ is diffuse. Hence, it suffices to prove that each $\Gamma$ invariant probability measure on $X$ is atomless.

    \medskip
    
    However, the space of $\Gamma$ invariant probability measures on $X$ is a Choquet simplex, with extreme points the ergodic measures $\text{Ext}(\text{Prob}_G(X))$ (See \cite{TC} 2.6). Therefore, for each $\Gamma$ invariant measure $\nu$ on $X$ there exists a measure $\eta$ on the space of ergodic measures $\text{Ext}(\text{Prob}_G(X))$ such that for all continuous functions $f$ on $X$ there holds:
    \begin{align*}
        \int_X f d\mu = \int_{\text{Ext}(\text{Prob}_G(X))} \int_X f(x) d\nu(x) d\eta(\nu).
    \end{align*}
    By an application of Urysohn's Lemma and the Dominated Convergence Theorem, one obtains the following for each closed subset $K$ of $X$:
    \begin{align*}
        \mu(K) = \int_{\text{Ext}(\text{Prob}_G(X))} \nu(K) d\eta(\nu).
    \end{align*}
    Now, recall that $\mu$ is a regular Borel measure, such that $\mu$ is atomless if and only if it assigns each singleton in $X$ mass $0$. Hence, applying the above, one has the following for each point $x_0$ in $X$
    \begin{align*}
        \mu(\{x_0\}) = \int_{\text{Ext}(\text{Prob}_G(X))} \nu(\{x_0\}) d\eta(\nu) = 0.
    \end{align*}
    Indeed, each Ergodic measure is atomless. Since $\mu$ was arbitrary, the proof is complete. 

\end{proof}

The condition that each ergodic measure be atomless may be viewed as a form of aperiodicity for the action. Indeed, suppose $\Gamma$ acts on $X$ by homeomorphisms and that $\mu$ is a Borel probability measure on $X$. Suppose additionally that $x_0$ is a point of $X$ such that $\mu(\{x_0\})>0$. By ergodicity of $\mu$, the orbit of $x_0$ has mass $1$. Hence, by subadditivity and invariance of the measure $\mu$ there holds:
\begin{align*}
    |O_\Gamma(x_0)| \cdot \mu(\{x_0\}) = 1,
\end{align*}
where $O_\Gamma(x_0)$ is the $\Gamma$ orbit of $x_0$. That is, an atom of $\mu$ gives rise to a finite $\Gamma$ orbit, and the size of the orbit may be computed using the measure of the point. Applying the contrapositive, if the action of $\Gamma$ on $X$ is free, then each ergodic measure is atomless. This gives rise to applications such as the following Corollary, which partially recovers a result of Li and Renault \cite{Xin}.

\medskip

\noindent \textbf{Corollary 12.2:} There is no decomposition:
\begin{align*}
    C_\lambda^*(\F_2) \cong C(X) \rtimes_{\alpha,r} \Gamma
\end{align*}
of the reduced group $C^*$-algebra of $\F_2$ as a crossed product arising from a free action of a countably infinite, discrete group $\Gamma$ by homeomorphisms on a compact, metrizable space $X$.
\begin{proof}
    Recall that the $1$-bounded entropy of the group von Neumann algebra $L(\F_2)$ is infinite (\cite{Voic2}, \cite{Jung}, and \cite{Reg}). Hence, $h(C_\lambda^*(\F_2)) = \infty$ by the variational principal. On the other hand, the above Corollary yields that the entropy of the reduced crossed product:
    \begin{align*}
        C(X) \rtimes_{\alpha,r} \Gamma
    \end{align*}
    vanishes.
\end{proof}

Observe also that the above Corollary yields a distinct computation for the entropy of irrational rotation algebras.

\medskip

Aside from direct crossed product decompositions, more flexible notions of normalizers for $C^*$-subalgebras are present in the literature. The above proof generalizes to this setting. Consider for instance the following Definition which appears in work of Renault \cite{Cartan}:

\medskip

\noindent \textbf{Definition 12.1 (\cite{Cartan} section 4):} Let $\mathscr{B}$ be a $C^*$-subalgebra of a $C^*$ algebra $\mathscr{A}$. Then the normalizer of $\mathscr{B}$, denoted $N_\mathscr{A}(\mathscr{B})$, is defined as follows:
\begin{align*}
    N_\mathscr{A}(\mathscr{B}) \defeq \{v \in \mathscr{A} \, | \, v^*\mathscr{B}v \subset \mathscr{B} \, \, \text{and} \, \, v \mathscr{B} v^* \subset \mathscr{B} \}. 
\end{align*}
Additionally, $\mathscr{B}$ is called regular provided its normalizer $N_\mathscr{A}(\mathscr{B})$ generates $\mathscr{A}$ as a $C^*$-algebra. 

\medskip

Once again taking as input an appropriate assumption about diffuse GNS representations, a $C^*$ version of the Theorems on normalizers present in work of Hayes (\cite{Reg} section 3) may be given. The author thanks David Jekel for offering a helpful outline of the following proof. 

\medskip

\noindent \textbf{Theorem 12.2:} Let $\mathscr{B}$ be a regular $C^*$-subalgebra of a $C^*$-algebra $\mathscr{A}$. Suppose additionally that for each trace $\tau$ on $\mathscr{A}$, the GNS von Neumann algebra $\pi_\tau(\mathscr{B})''$ is diffuse. Then:
\begin{align*}
    h(\mathscr{A}) \leq h(\mathscr{B}).
\end{align*}
\begin{proof}
    Fix a trace $\tau$ on $\mathscr{A}$. Fix $v$ belonging to $N_\mathscr{A}(\mathscr{B})$ and consider $x \defeq \pi_\tau(v)$. Denote $Q \defeq \pi_\tau(\mathscr{B})''$ and $M \defeq \pi_\tau(\mathscr{A})''$. It will be demonstrated that $x$ belongs to the anti-coarse space (see Hayes, Jekel, Kunnawalkam Elayavalli \cite{PT} section 2 or \cite{Reg} Proposition 3.2 for more details):
    \begin{align*}
        L^2_\dagger(Q \leq M) \defeq \bigcap_{T \in \text{Hom}_{Q-Q}(L^2(M),L^2(Q) \otimes L^2(Q))} \ker(T).
    \end{align*}
    For this, set $T$ to be a $Q-Q$ bimodular map between $L^2(M)$ and the coarse bimodule $L^2(Q) \otimes L^2(Q)$. Observe that since $\mathscr{B}$ is unital and $v$ is in the normalizer, $xx^*$ and $x^*x$ both belong to $Q$. Writing the polar decomposition $x = h|x|$, it follows that both the source and range projections $q$ and $p$ respectively of $h$ belong to $Q$.

    \medskip
    
    Since $Q$ is diffuse, there exists a sequence $u_n$ of unitaries in $pQp$ which converge to $0$ in the weak operator topology. Similarly, pick out unitaries $w_n \defeq hu_nh^*$ in $qQq$ which converge to $0$ in the weak operator topology. Now, by bimodularity:
    \begin{align*}
        \|T(h)\|_2^2 = \langle T(w_n^*w_n h),T(h))\rangle = \langle w_n^*T(h)u_n,T(h) \rangle \rightarrow 0.
    \end{align*}
    Indeed, both $\{u_n\}$ and $\{w_n\}$ converge to $0$ in the weak operator topology. Hence, $T(h)=0$. Applying bimodularity once more to the polar decomposition (note that $|x|$ belongs to $Q$) one obtains that $T(x)=0$. Hence, $x$ belongs to the anti-coarse space. By the regularity assumption,
    \begin{align*}
        \pi_\tau(N_\mathscr{A}(\mathscr{B}))'' = M.
    \end{align*}
    Hence, $W^*(L^2_\dagger(Q \leq M)) = M$. Therefore, applying Hayes \cite{Reg} Theorem 1.5 yields:
    \begin{align*}
        h(M) \leq h(Q).
    \end{align*}
    Recalling the definitions of $M$ and $Q$ and applying the variational principle completes the proof. 
\end{proof}

\section{Tensor Products}

In the regularity paper of Hayes, one of the chief applications of $1$-bounded entropy is the proof of primeness for positive entropy II$_1$ factors \cite{Reg}. In this section, the analogous primeness results for $C^*$-algebras are obtained. 

\medskip

\noindent \textbf{Theorem 13.1:} Suppose that $\mathscr{A}$ is a unital $C^*$ algebra such that $\pi_\tau(\mathscr{A})''$ is diffuse for each trace $\tau$ belonging to $\mathscr{T}(\mathscr{A})$. Similarly, suppose $\mathscr{B}$ is a unital $C^*$ algebra such that $\pi_\tau(\mathscr{B})''$ is diffuse for each trace $\tau$ belonging to $\mathscr{T}(\mathscr{B})$. Then 
\begin{align*}
    h(\mathscr{A} \otimes_\text{min} \mathscr{B}) \leq 0.
\end{align*}
\begin{proof}
    Fix a trace $\tau$ on the minimal tensor product $\mathscr{A} \otimes_\text{min} \mathscr{B}$. By the unital assumption on both $\mathscr{A}$ and $\mathscr{B}$, there are well-defined factor embeddings:
    \begin{align*}
       & \mathscr{A} \hookrightarrow \mathscr{A} \otimes_\text{min} \mathscr{B}  \\
       & \mathscr{B} \hookrightarrow \mathscr{A} \otimes_\text{min} \mathscr{B}.
    \end{align*}
    Therefore, it makes sense to consider the restrictions $\tau|_\mathscr{A}$ and $\tau_\mathscr{B}$ of the trace. Observe that the GNS von Neuman algebra $\pi_\tau(\mathscr{A} \otimes_\text{min} \mathscr{B})''$ is generated by diffuse, commuting algebras isomorphic to $\pi_{\tau|_\mathscr{A}}(\mathscr{A})''$ and $\pi_{\tau|_\mathscr{B}}(\mathscr{B})''$ respectively. This is due to the uniqueness of GNS representations. Therefore, since the tensor product of diffuse tracial von Neumann algebras has $1$-bounded entropy at most $0$ \cite{Reg}, applying the variational principle yields the result. 
\end{proof}

This yields indecomposability Corollaries for $C^*$ algebras such at $C_\lambda^*(\F_2)$, which has infinite entropy. Recall again that some $C^*$-algebras obeying the "diffuseness" hypothesis in the Theorem include nowhere scattered $C^*$ algebras as well as the reduced group $C^*$-algebras of non-amenable groups. The example in the following section on quotients indicates the necessity that both tensorands satisfy such a condition. Indeed, it exhibits a minimal tensor product of infinite entropy.

\section{Joins and Free Products}

In this section, the generalization of the join Lemma \cite{Reg} A.12 is achieved along with some free product indecomposability results. 

\medskip

\noindent \textbf{Theorem 14.1:} Suppose that $\mathscr{A}$ and $\mathscr{B}$ are $C^*$-algebras represented on a common Hilbert space $\mathscr{H}$. Consider the $C^*$ algebra $\mathscr{D}$ generated by $\mathscr{A}$ and $\mathscr{B}$, and suppose that the intersection $\mathscr{C} \defeq \mathscr{A} \cap \mathscr{B}$ is such that for each trace $\tau$ on $\mathscr{D}$, the GNS representation $\pi_{\tau|_\mathscr{C}}(\mathscr{C})''$ is diffuse. Then:
\begin{align*}
    h(\mathscr{D}) \leq h(\mathscr{A})+h(\mathscr{B}).
\end{align*}
\begin{proof}
    Fix a trace $\tau$ on $\mathscr{D}$, and observe that $\pi_\tau(\mathscr{D})''$ is the join of the tracial von Neumann algebras $\pi_{\tau|_\mathscr{A}}(\mathscr{A})''$ and $\pi_{\tau|_\mathscr{B}}(\mathscr{B})''$ which have diffuse intersection. Indeed, $\pi_{\tau|_\mathscr{C}}(\mathscr{C})''$ is diffuse. Hence, applying the Join Lemma \cite{Reg} A.12 yields:
    \begin{align*}
        h(\pi_\tau(\mathscr{D})'') \leq h(\pi_{\tau|_\mathscr{A}}(\mathscr{A})'')+h(\pi_{\tau|_\mathscr{B}}(\mathscr{B})'').
    \end{align*}
    Now, apply the variational principle to obtain the following:
    \begin{align*}
        h(\mathscr{D})= \sup_{\tau \in \mathscr{T}(\mathscr{D})} h(\pi_\tau(\mathscr{D})'') \leq \sup_{\tau \in \mathscr{T}(\mathscr{D})} [h(\pi_{\tau|_\mathscr{A}}(\mathscr{A})'')+h(\pi_{\tau|_\mathscr{B}}(\mathscr{B})'') ] \\\leq \sup_{\tau \in \mathscr{T}(\mathscr{A})}h(\pi_\tau(\mathscr{A})'') + \sup_{\tau \in \mathscr{T}(\mathscr{B})}h(\pi_\tau(\mathscr{B})'') = h(\mathscr{A})+h(\mathscr{B}).
    \end{align*}
\end{proof}

On the other hand, full free products typically yield infinite entropy. 

\medskip

\noindent \textbf{Proposition 14.1:} Suppose that $\mathscr{A}$ and $\mathscr{B}$ are $C^*$ algebras equipped with traces such that at least one of $\dim \mathscr{A}$ or $\dim \mathscr{B}$ exceeds $2$. Then 
\begin{align*}
    h(\mathscr{A}*\mathscr{B})=\infty,
\end{align*}
where the free product is the universal free product. 
\begin{proof}
    Fix a trace $\tau_\mathscr{A}$ on $\mathscr{A}$ and a trace $\tau_\mathscr{B}$ on $\mathscr{B}$. Consider the free product von Neumann algebra $W^*(\mathscr{A},\tau_\mathscr{A})*W^*(\mathscr{B},\tau_\mathscr{B})$. By the universal property of the full free product, this is the von Neumann algebra generated by a GNS representation of $\mathscr{A}*\mathscr{B}$ with respect to a trace. This algebra has infinite $1$-bounded entropy by results of Brown, Dykema, and Jung \cite{FPL}.
\end{proof}

In the exceptional case in which each algebra has dimension $2$, the computation is as follows.

\medskip

\noindent \textbf{Example 14.1:} Consider the free product $\mathscr{A} * \mathscr{B}$ in the case in which both $\mathscr{A}$ and $\mathscr{B}$ coincide with $\C^2$. Observe:
\begin{align*}
    \C^2 * \C^2 = C^*(\Z/(2) * \Z/(2)).
\end{align*}
However, $\Z/(2)*\Z/(2)$ is amenable, such that $C^*(\Z/(2) * \Z/(2)) = C_\lambda^*(\Z/(2)*\Z/(2))$, the reduced group $C^*$ algebra. Since $\Z/(2)*\Z/(2)$ is amenable, this algebra is nuclear, and so has entropy $0$ by the above.

\section{Behavior Under Quotients}
Recall that the only ultraweakly closed ideals within a von Neumann algebra are corners by central projections, I.e. direct summands. Hence, the problem of understanding the behavior of $1$-bounded entropy under quotients is uninteresting in the setting of von Neumann algebras. The situation changes drastically in the setting of $C^*$ algebras, where the ideal structure is rich. The following Proposition has already been used tacitly above in a couple of examples.

\medskip

\noindent \textbf{Proposition 15.1:} Suppose that $I \subset \mathscr{A}$ is a norm closed $2$-sided ideal in a $C^*$ algebra $\mathscr{A}$. Then $h(\mathscr{A}/I) \leq h(\mathscr{A})$.
\begin{proof}
    Denote the quotient map $q: \mathscr{A} \rightarrow \mathscr{A}/I$. Each trace $\tau$ on the quotient $\mathscr{A}/I$ defines a trace on $\mathscr{A}$ by $\tau \circ q$ with:
    \begin{align*}
        \pi_{\tau \circ q}(\mathscr{A}) = \pi_\tau(\mathscr{A}/I).
    \end{align*}
    Now apply the variational principle. 
\end{proof}

This behavior was key to the above example regarding $C^*(\F_2)$. It also arises in the following example, which again reiterates the importance of the above notion of "diffuse" for $C^*$-algebras.

\medskip

\noindent \textbf{Example 15.1:} The reduced group $C^*$-algebra $C_\lambda^*(\F_2 \times \Z) \cong C_\lambda^*(\F_2) \otimes_\text{min} C(S^1)$ has infinite $1$-bounded entropy. 
\begin{proof}
    The group homomorphism $\F_2 \times \Z \rightarrow \F_2$ has amenable kernel and so extends to a surjective homomorphism 
    \begin{align*}
        C_\lambda^*(\F_2 \times \Z) \rightarrow C_\lambda^*(\F_2)
    \end{align*}
    of reduced group $C^*$ algebras. Recall now that $h(C_\lambda^*(\F_2))=\infty$ by the variational principle and the fact that the free group factor $L(\F_2)$ has infinite $1$-bounded entropy. Now apply the proposition on quotients. 
\end{proof}

\medskip

\noindent \textbf{Remark:} While quotients are under consideration, a little remark about the absence of microstates is pertinent. Recall that a von Neumann algebra has $1$-bounded entropy $-\infty$ provided it is not Connes embeddable. One might ask what it means for a $C^*$ algebra to have topological entropy $-\infty$. It turns out that this occurs if and only if the algebra has no quotient which is matricially finite. 

\section{Property (T) Groups}

In this section, the object of interest is the full group $C^*$ algebra of a property (T) group. This will give an example of a $C^*$-algebra for which it may be proved that its entropy if finite, but for which it remains unknown whether its entropy is $0$. The distinction between finite entropy and entropy $0$ for tracial von Neumann algebras is open.

\medskip

The treatment is analogous to Hayes, Jekel, and Kunnawalkam Elayavalli \cite{PropT}. Recall that if a countable, discrete group $G$ has property (T), then there exists a finite subset $F \subset G$ and a constant $\gamma>0$ (referred to as the Kazhdan constant) for which:
\begin{align*}
    \|\xi - P\xi\| \leq \frac{1}{\gamma} \sum_{g \in F} \|\pi(g)\xi-\xi\|
\end{align*}
for all unitary representations $\pi$ of $G$ on Hilbert space $\mathscr{H}$ and $\xi$ belonging to $\mathscr{H}$. Here, $P$ is the projection onto the space of $G$-invariant vectors within $\mathscr{H}$. The finite set $F$, if it is ordered, is referred to as a Kazhdan tuple for $G$. 

\medskip

In analogy to this, a tuple $(x_1,x_2,\dots,x_n)$ of self-adjoint operators within a tracial von Neumann algebra $(M,\tau)$ is referred to as a Kazhdan tuple provided there exists a constant $\gamma>0$ (still referred to as a Kazhdan constant) such that for every Hilbert $M-M$ bimodule $\mathscr{H}$ and vector $\xi$ within $\mathscr{H}$ there holds:
\begin{align*}
    \|\xi-P\xi\| \leq \frac{1}{\gamma} \|x\xi - \xi x\|_{\mathscr{H}^{\oplus n}}.
\end{align*}
Here, once again, $P$ denotes the projection onto the space of central vectors in $\mathscr{H}$. Moreover, 
\begin{align*}
    & x\xi = (x_1\xi,x_2 \xi, \dots, x_n \xi) \\
    & \xi x = (\xi x_1, \xi x_2, \dots, \xi x_n).
\end{align*}

The following Proposition from \cite{PropT} will be of use.

\medskip

\noindent \textbf{Proposition 16.1 (\cite{PropT} Lemma 2.7):} Suppose that $G$ is a property (T) group. Fix a Kazhdan tuple $g_1,g_2,\dots,g_r$ for $G$ with constant $\gamma>0$. Then, for each unitary representation $\pi$ of $G$ on a Hilbert space $\mathscr{H}$, one has that 
\begin{align*}
    (\text{Re} \, \pi(g_1),\text{Re} \, \pi(g_2), \dots, \text{Re} \, \pi(g_r), \text{Im} \, \pi(g_1),\text{Im} \, \pi(g_2), \dots, \text{Im} \, \pi(g_r))
\end{align*}
is a Kazhdan tuple for $W^*(\pi(G))$ with the same Kazhdan constant $\gamma>0$. 

\medskip

Henceforth, for $G$ a property (T) group, let $x$ be the self-adjoint tuple:
\begin{align*}
    (\text{Re} \, g_1,\text{Re} \, g_2, \dots, \text{Re} \, g_r, \text{Im} \, g_1,\text{Im} \, g_2, \dots, \text{Im} \, g_r)
\end{align*}
within $C^*(G)$ arising from a Kazhdan tuple for $G$ with constant $\gamma>0$. Observe that if $\pi$ is a unitary representation of $G$, then $\pi(x)$ is a Kazhdan tuple of length $2r$ for $W^*(\pi(G))$ with constant $\gamma>0$ from proposition 16.1. Moreover, denote its law $\ell_{\pi(x)}$. By the universal property of $C^*(G)$, one has that the set $L_x$ consists of the laws $\ell_{\pi(x)}$ as $\pi$ ranges over the unitary dual of $G$. 

\medskip

Now, the main result of the section may be demonstrated.

\medskip

\noindent \textbf{Theorem 16.1:} Let $G$ be a countable, discrete group with property (T). Then $h(C^*(G))< \infty$.
\begin{proof}
    Adopt the notation of the foregoing discussion. Let $R>0$ be such that each member $x_i$ of the $2r$-tuple $x$ satisfies $\|x_i\| < R$. Since a $*$-homomorphism between $C^*$ algebras is contractive, the same norm estimate holds for the tuples $\pi(x)$ for any unitary representation $\pi$ of $G$ (notation is abused here as $\pi$ and its unique extension to $C^*(G)$ are both denoted $\pi$). Fix $0<\epsilon<\min(1,\gamma)$ and define for each compact, convex subset $K \subset \Sigma_{R,2r}$ the quantity:
    \begin{align*}
        h_{\epsilon,R}(K) \defeq \inf_{\mathcal{O} \supset K} \limsup_{N \rightarrow \infty} \frac{1}{N^2}\log K_\epsilon^\text{orb}(\Gamma_R^{(N)}(\mathcal{O})).
    \end{align*}
    Here, the notation indicating the finite subset of the index set is suppressed. Indeed, the index set is finite of size $2r$, so all packing is done with respect to the whole index set. Similarly, for each unitary representation $\pi$ of $G$, denote:
    \begin{align*}
        h_{\epsilon,R}(\pi(x)) \defeq h_{\epsilon,R}(\{\ell_{\pi(x)}\}) = h_{\epsilon,R}(W^*(\pi(x))).
    \end{align*}
    Examining the compactness argument in the proof of the variational principle, for fixed $\epsilon>0$ there exists a law $\ell$ belonging to $L_x$ such that:
    \begin{align*}
        h_{\epsilon/2,R}(\ell) = h_{\epsilon,R}(L_x).
    \end{align*}
    Indeed, the reason why the supremum in the general variational principle is not a priori realized is due to the final supremum over $\epsilon$ (and, in general, over finite subsets) alone. Recall that this law $\ell$ is the law $\ell_{\pi(x)}$ for some unitary representation $\pi$ of $G$. Now, let $\rho$ be an arbitrary unitary representation of $G$. Applying the estimate in the proof of \cite{PropT} 1.1, there exists $C>0$ such that:
    \begin{align*}
        h(W^*(\rho(x))) \leq h_{\epsilon,R}(\rho(x)) + \frac{12(2r+1)}{\gamma} \sum_{j=0}^\infty \epsilon^{2^j} \bigg(\log(CR(2r)^{1/2})+2^{j+1} \log(1/\epsilon) \bigg). 
    \end{align*}
    By design, $h_{\epsilon,R}(\rho(x)) \leq h_{\epsilon,R}(L_x) = h_{\epsilon/2,R}(\pi(x))$. Now, by \cite{PropT} Theorem 1.1: 
    \begin{align*}
        h_{\epsilon/2,R}(\pi(x)) \leq h(W^*(\pi(x))) < \infty.
    \end{align*}
    Moreover, the above sum converges and does not depend on $\rho$. Indeed, it depends only on the Kazhdan constant $\gamma>0$ and the length $2r$. Therefore, a \textit{uniform} upper bound on $h(W^*(\rho(x)))$ over all unitary representations $\rho$ of $G$ has been established. Hence, applying the variational principle yields $h(C^*(G)) < \infty$. 
    
\end{proof}

\noindent \textbf{Remark:} Consequently, by Theorem 8.1, one also has $h^\text{top}(C^*(G)) < \infty$. The topological entropy on $C^*(G)$ seems more natural given the universal property of $C^*(G)$.

\printbibliography

\end{document}